\documentclass[basic]{custom}
 \usepackage{times, mathptm}
 \usepackage{latexsym, amsmath, amssymb, bm, mathtools} 
 \usepackage{hyperref}
 \definecolor{darkblue}{rgb}{0,0,0}
 \hypersetup{colorlinks=true, linkcolor=darkblue,
   urlcolor=darkblue, citecolor=darkblue}
 \usepackage{multirow} \usepackage[authoryear]{natbib}
 \usepackage{float}
 \usepackage{wasysym}

 \usepackage{fancyvrb}
 \usepackage[scale=0.9293]{zlmtt}
 \def\CODE{\begin{Verbatim}[fontsize=\footnotesize,fontseries=b]}

 \def\NI{\noindent}
 \def\TS{\textstyle} \def\d{\textrm{d}} \def\i{\textrm{i}} 
 \def\VEC#1{{\bm{#1}}} \def\MAT#1{{\bm{#1}}}
 \def\IE{{\it i.e.}} \def\EG{{\it e.g.}}  
 \def\EEG{{\it E.g.}}
   
 \def\THIRD{{$3^{\textrm{rd}}$}}  
 \def\DOTS{{...}} \def\DS{\displaystyle}
 \def\Z{\mathbf{Z}} \def\R{\mathbf{R}}

 \def\STRUT{\vphantom{\vert_\vert^\vert}}

 \def\SEP{,\mkern1.5mu} \def\DWS{\mkern3mu}
 \def\AIJ{{a_{\mkern1.15mu i \mkern-0.85mu j}}}
 \def\AIN#1{a_{\mkern1mu i \mkern1mu #1}}
 \def\BI{{b_{\mkern0.8mu i}}}
 \def\BJ{{b_{\mkern-1.25mu j}}}
 \def\BCN#1{\VEC{b} \mkern0.75mu \VEC{c}^{#1}}
 \def\CI{{c_{\mkern0.5mu i}}}
 \def\CIN{{c_{\mkern0.5mu i}^{\mkern0.5mu n}}}
 \def\CJ{{c_{\mkern-1.25mu j}}}

 \def\DN{\mathcal{D}_{\mkern1mu n}}

 \def\DO{{\VEC{d}_{\mkern-0.5mu{}1}}}
 \def\DOA{{\VEC{d}_{\mkern-0.5mu{}1} \mkern1mu \MAT{A}}}
 
 \def\FI{{\VEC{F}_{\mkern-1.5mu i}}}
 \def\FJ{{\VEC{F}_{\mkern-4mu j}}}

 \def\DT{\d \mkern1mu t}
 \def\X12S#1{{\VEC{X}_{\mkern-1mu{#1}}}}
 \def\APHI{{\MAT{A} \mkern0.25mu \VEC{\Phi}}}
 \def\BPHI{{\VEC{b} \mkern1.75mu \VEC{\Phi}}}
 \def\BRAT{( \mkern0.25mu t \mkern-0.5mu )}
 \def\BRA1{(\VEC{1} \mkern-0.75mu )}
 \def\BRAA1{(\MAT{A} \VEC{1} \mkern-0.75mu )}
 \def\BRAP1{(\VEC{p}_1 \mkern-0.5mu )}
 \def\BRAAP1{(\MAT{A} \VEC{p}_1 \mkern-0.5mu )}
 \def\BRAAQ1{(\MAT{A} \mkern0.5mu \VEC{q}_1 \mkern-0.5mu )}
 \def\INT01{$[ \mkern1mu 0, 1 \mkern0.5mu ]$}
 \def\TTM{\textrm{t}_{\mkern1mu m}}
 \def\TTN{\textrm{t}_{\mkern1mu n}}
 \def\ORDT{\vert \mkern1mu \textrm{t} \mkern1mu \vert}
 \def\TBPN{[ \mkern1.75mu \bullet^n ]}
 \def\ATT{a(\textrm{t}) \mkern1mu \textrm{t}}
 \def\TPRIME{\textrm{t}\mkern0.75mu{}'}
 \def\BCDQ{CNC}
 \def\QED{$\Box$}
 \def\PONE{{\small\textsf{1}}}
 \def\PTWO{{\small\textsf{2}}}
 \def\PTHREE{{\small\textsf{3}}}
 \def\PFOUR{{\small\textsf{4}}}
 \def\PFIVE{{\small\textsf{5}}}
 \def\PSIX{{\small\textsf{6}}}
 \def\PSEVEN{{\small\textsf{7}}}
 \def\PEIGHT{{\small\textsf{8}}}
 
 \def\SLASH{$\mkern2mu$/$\mkern2mu$}
 \def\FIELD{\pmb{\mathbf{Q}}(\alpha \SEP \beta)}
 \def\SUBD#1{D_{\mkern1mu #1}}

\begin{document}

 \title{On~Runge--Kutta~methods~of~order~10}

 \author*{\hskip54pt\fnm{Misha} \sur{Stepanov}
\qquad {\footnotesize\tt stepanov@arizona.edu}}

 \affil{\orgdiv{Department of Mathematics \,\;and\;\, Program in Applied 
Mathematics}, \orgname{University of Arizona}, 
\orgaddress{\city{Tucson}, \state{AZ 85721}, 
\country{US$\mkern0.5mu$A}}}

\abstract{\hskip10pt{}A family of explicit $15$-stage Runge--Kutta 
methods of order $10$ is derived.~~~~~}

\keywords{minimal number of stages, explicit Runge--Kutta methods}

\pacs[MSC Classification]{65L05, 65L06}

\maketitle

\hskip\parindent{}Runge--Kutta methods (see, \EG, \citep[s.~23 and 
ch.~3]{But16}, \citep[ch.~II]{HNW93}, \citep[ch.~4]{AsPe98}, 
\citep[ch.~3]{Ise08}) are widely and successfully used to solve ordinary 
differential equations numerically for over a century \citep[]{BuWa96}. 
Being applied to a system $\d \mkern0.5mu \VEC{x} / \DT = \VEC{f}(t \SEP 
\VEC{x})$, in order to propagate by the step size $h$ and update the 
position, $\VEC{x}\BRAT \mapsto \VEC{\tilde{x}}(t + h)$, where 
$\VEC{\tilde{x}}(t + h)$ is a numerical approximation to the exact 
solution $\VEC{x}(t + h)$, an $s$-stage Runge--Kutta method (which is 
determined by the coefficients $\AIJ$, weights $\BJ$, and nodes 
$c_{\mkern0.5mu i}$) would form the following system of equations for 
$\X12S{1}$, $\X12S{2}$, \DOTS, $\X12S{s}$:

\vspace{-12pt}

 \begin{gather*}
   \X12S{i} = \VEC{x}\BRAT + h \mkern-1.5mu \sum_{j = 1}^{s} \AIJ 
\mkern2mu \FJ , \qquad \FI = \VEC{f} \bigl( t + \CI \mkern1mu h \SEP 
\X12S{i} \bigr), \qquad i = 1, \mkern1.5mu 2, \DOTS, \mkern1.5mu s
 \end{gather*}

\vspace{-6pt}

\NI solve it, and then compute $\VEC{\tilde{x}}(t + h) = \VEC{x}\BRAT + 
h \sum_{\mkern-1mu j = 1}^s \BJ \mkern2mu \FJ$. In the limit $h \to 0$ 
all the vectors $\FI$, where $1 \le i \le s$, are the same, so it is 
natural and will be assumed that $\sum_{\mkern-1mu j = 1}^{s} \AIJ = 
\CI$ for all $i$.\footnote{~See \citep[eq.~(3.8)]{Oli75} for an example 
of a $2$-stage Runge--Kutta method of order $2$ that violates this 
assumption.}

A method is said to be of \emph{order} [at least] $p$ if for sufficiently 
smooth r.h.s. function $\VEC{f}$ the local truncation error behavior is 
$\Vert \VEC{x}(t + h) - \VEC{\tilde{x}}(t + h) \Vert = 
\mathcal{O}(\smash{h^{\mkern1mu p + 1}})$ as $h \to 0$. It is often 
desirable to use a method of higher order, as that allows to obtain a 
solution with a certain level of accuracy with a smaller number of steps. 
For an $s$-stage Runge--Kutta method the maximal possible order is $p = 2 
s$, achieved by Gauss--Legendre methods \citep{But64a}.

A Runge--Kutta method is called \emph{explicit} if $\AIJ = 0$ whenever 
$j \ge i$. Then $c_1 = 0$, $\X12S{1} = \VEC{x}\BRAT$, $\VEC{F}_1 = 
\VEC{f} \bigl( t \SEP \VEC{x}\BRAT \bigr)$, $a_{21} = c_2$, and 
$\X12S{2}$, $\VEC{F}_2$, $\X12S{3}$, $\VEC{F}_3$, \DOTS, $\X12S{s}$, 
$\VEC{F}_s$ could be computed in sequence by direct computation. \EEG, 
at the moment of finding $\X12S{3}$ the vectors $\X12S{1}$, $\VEC{F}_1$, 
$\X12S{2}$, $\VEC{F}_2$ are already computed, and

\vspace{-16pt}

 \begin{gather*}
    \X12S{3} \,=\, \VEC{x}\BRAT + h \mkern1mu a_{31} \mkern2mu \overbrace{ \mkern-2mu \VEC{f} \bigl( t \SEP 
\mkern-3mu \underbrace{\mkern-1mu \VEC{x}\BRAT \mkern-3mu }_{\X12S{1}} \mkern-1mu \bigr) \mkern-2mu }^{\VEC{F}_1} + h \mkern1mu 
a_{32} \mkern2mu \overbrace{ \mkern-2mu \VEC{f} 
\mbox{\raisebox{-1pt}{\scalebox{1.2}{$\bigl($}}} \mkern1mu
 t + c_2 h \SEP \underbrace{\VEC{x}\BRAT + h \mkern1mu c_2 
\mkern1mu \VEC{f} \bigl( t \SEP \VEC{x}\BRAT \bigr) \mkern-3mu}_{\X12S{2}} 
 \mbox{\raisebox{-1pt}{\scalebox{1.2}{$\bigr)$}}} \mkern-2mu }^{\VEC{F}_2}
 \end{gather*}

\vspace{-7pt}

\NI{}Determining the minimal number of stages $s_{\textrm{min}}(p)$ for 
which there exists an explicit Runge--Kutta method of order $p$ is a 
complicated problem, which is currently solved for $p \le 8$: 
$s_{\textrm{min}}(\langle 1 \SEP 2 \SEP 3 \SEP 4 \SEP 5 \SEP 6 \SEP 7 
\SEP 8 \mkern1mu \rangle) = \langle 1 \SEP 2 \SEP 3 \SEP 4 \SEP 6 \SEP 7 
\SEP 9 \SEP 11 \mkern0.5mu \rangle$, with the lower bound 
$s_{\textrm{min}}(p) \ge p + 3$ for $p > 8$ \citep{But85}.

There are known explicit methods of order $10$ with $18$ stages 
\citep{Cur75}; with $17$ stages: \citep{Hai78}, following its structure 
\citep{Ono03}, and \citep{Fea07} with performance traded off for the 
presence of an embedded method of order~$8$; and with $16$ stages 
\citep{Zha24} (although there is no yet a rigorous proof that the method 
is indeed of order $10$, the numerical evidence is overwhelming).

The aim of this work is to construct an explicit $15$-stage Runge--Kutta 
method of order $10$. Order conditions are stated in 
Section~\ref{order_conditions}. Order conditions of two types, Q- and 
D-types, are considered in Section~\ref{QDOC}, while in 
Sections~\ref{dual_method} and \ref{heuristics} these are compared and 
contrasted. A $7$-dimensional family of explicit $15$-stage Runge--Kutta 
methods of order $10$ is derived in Section~\ref{family}. Some previously 
known methods of order $10$ are compared to a selected new one in 
Section~\ref{comparison}.

\section{Order conditions} \label{order_conditions}

\hskip\parindent{}The element-wise product of tensors $\VEC{x}$ and 
$\VEC{y}$ of the same size will be denoted as $\VEC{x} \mkern0.5mu . 
\mkern1mu \VEC{y}$, \EG, in case of vectors $(\VEC{x} \mkern0.5mu . 
\mkern1mu \VEC{y})_i = x_i \mkern1mu y_i$. The element-wise product of 
$n$ copies of a column vector $\VEC{x}$ will be written as 
$\VEC{x}^{\mkern1mu n}$. Let $\VEC{1}$ be the $s$-dimensional column 
vector with all components being equal to $1$; $\MAT{A} = \bigl[ \AIJ 
\bigr]$ be the $s \times s$ matrix with $\AIJ$ as its matrix element in 
the $i^{\mkern1mu \textrm{th}}$ row and $j^{\mkern1mu \textrm{th}}$ 
column; $\VEC{b} = \bigl[ \BJ \bigr]$ be the weights row vector; and 
$\VEC{c} = \bigl[ \CI \bigr]$ be the nodes column vector.

Given rooted trees $\textrm{t}_1$, $\textrm{t}_2$, \DOTS, $\TTN$, a new 
tree $[ \mkern1mu \textrm{t}_1 \mkern2mu \textrm{t}_2 \mkern2mu \DOTS 
\mkern2mu \TTN \mkern1mu ]$ is obtained by connecting with $n$ edges 
their roots to a new vertex, the latter becomes a new root 
\citep[s.~301]{But16}, \citep[p.~152]{HNW93}, \citep[p.~44]{But21}, 
\citep[p.~53]{HLW06}. Consider a vector function $\VEC{\Phi} : \mathrm{T} 
\to \pmb{\mathbf{R}}^s$ on the set of rooted trees that is recursively 
defined as $\VEC{\Phi}( \bullet ) = \VEC{1}$ and $\VEC{\Phi} \bigl( [ 
\mkern1mu \textrm{t}_1 \mkern2mu \textrm{t}_2 \mkern2mu \DOTS \mkern2mu 
\TTN \mkern1mu ] \bigr) = \prod_{\mkern1mu m = 1}^{\mkern1mu n} 
\APHI(\TTM)$, where the product of vectors is taken element-wise. This 
function coincides with \emph{derivative weights} 
\citep[def.~312A]{But16}, \citep[pp.~148 and 151]{HNW93}, it is closely 
related to \emph{internal} or \emph{stage weights} $\APHI(\textrm{t})$ 
\citep[p.~125]{But21} and \emph{elementary weights} $\BPHI(\textrm{t})$ 
\citep[p.~55]{HLW06}.

A Runge-Kutta method $(\MAT{A} \SEP \VEC{b} \SEP \VEC{c})$ is of order 
[at least] $p$ if and only if for any rooted tree $\textrm{t}$, with 
$\ORDT \le p$, one has $\BPHI(\textrm{t}) = 1 / \textrm{t}!$ 
\citep[s.~315]{But16}, \citep[p.~153]{HNW93}, \citep[p.~127]{But21}, 
\citep[p.~56]{HLW06}. Here $\ORDT$ is the order of tree $\textrm{t}$, 
\IE, the number of vertices in $\textrm{t}$. The factorial $\textrm{t}!$ 
is recursively defined as $\bullet ! = 1$ and if $\mathrm{t} = [ 
\mkern1mu \textrm{t}_1 \mkern2mu \textrm{t}_2 \mkern2mu \DOTS \mkern2mu 
\TTN \mkern1mu ]$, then $\mathrm{t}! = \ORDT \mkern1mu \prod_{\mkern1mu m 
= 1}^{\mkern1mu n} \bigl( \TTM \bigr) \mkern-0.5mu ! \mkern1mu$.

\section{Q- and D-type order conditions} \label{QDOC}

\hskip\parindent{}For any rooted tree $\textrm{t}$, let 
$\VEC{Q}(\textrm{t}) = \APHI (\textrm{t}) - \VEC{c}^{\ORDT} \mkern-1mu / 
\textrm{t}!$. For a rooted tree $\textrm{t} = \TBPN$ of height at most 
$1$, the vector $\VEC{Q} \bigl( \TBPN \bigr) = \VEC{q}_n = \MAT{A} 
\VEC{c}^{\mkern0.5mu n} - \frac{1}{n + 1} \VEC{c}^{\mkern0.5mu n + 1}$ is 
a standard subquadrature vector \cite[p.~1124]{VeZe95}, 
\cite[p.~558]{Ver14}. It is convenient to define $[ \mkern1.75mu 
\bullet^0 \mkern0.5mu ] = [ \varnothing \mkern0.25mu ] = \bullet$. The 
assumed condition $\MAT{A} \VEC{1} = \VEC{c}$ coincides with 
$\VEC{Q}(\bullet) = \VEC{q}_0 = \VEC{0}$.

Consider a vector space of all finite linear combinations of trees 
$\mathcal{T} = \oplus_{\mkern1.5mu \textrm{t} \in \textrm{T}} \mkern1mu 
\R$. Such a linear combination could be written as 
${\footnotesize\textsf{T}} = \sum_{\mkern1.5mu \textrm{t} \in 
\textrm{T}} \ATT$, with the weights $a(\textrm{t})$ being non-zero only 
for finitely many trees $\textrm{t}$. Let $\mathcal{Q} : \textrm{T} \to 
\mathcal{T}$ be a mapping defined as $\mathcal{Q}(\textrm{t}) = [ 
\mkern2.25mu \textrm{t} \mkern2.25mu ] - \frac{1}{\textrm{t}!} [ 
\mkern1.75mu \bullet^{\ORDT} \mkern0.5mu ]$. One has 
$\mathcal{Q}(\bullet) = 0$. All three mappings $\VEC{\Phi}$, $\VEC{Q}$, 
and $\mathcal{Q}$ could be extended to linear combinations of trees by 
linearity: $\langle \VEC{\Phi} \SEP \VEC{Q} \SEP \mathcal{Q} \mkern1mu 
\rangle \bigl( \sum_{\textrm{t} \in \textrm{T}} \ATT \bigr) = 
\sum_{\textrm{t} \in \textrm{T}} a(\textrm{t}) \mkern1mu \langle 
\VEC{\Phi} \SEP \VEC{Q} \SEP \mathcal{Q} \mkern1mu \rangle 
(\textrm{t})$. One has $\VEC{Q}({\footnotesize\textsf{T}}) = \VEC{\Phi} 
\bigl( \mkern-1.25mu \mathcal{Q}({\footnotesize\textsf{T}}) \bigr)$ for 
any ${\footnotesize\textsf{T}} \in \mathcal{T}$.

The order conditions $\BPHI(\textrm{t}) = 1 / \textrm{t}!$ whenever 
$\ORDT \le p$ could be rewritten as quadrature conditions $\VEC{b} 
\VEC{c}^{\mkern0.5mu n} = \frac{1}{n + 1}$ for $0 \le n < p$ and order 
conditions of \emph{Q-type} $\VEC{b} \mkern1mu (\VEC{Q}( \mkern0.5mu 
\textrm{t}_1 \mkern-0.5mu ) \mkern0.5mu . \mkern1mu \VEC{Q}( \mkern0.5mu 
\textrm{t}_2 \mkern-0.5mu ) \mkern0.5mu . 
\mkern-2mu\cdot\mkern-2.5mu\cdot\mkern-2.5mu\cdot\mkern-2mu . \mkern1mu 
\VEC{Q}( \mkern0.5mu \textrm{t}_k \mkern-0.5mu ) \mkern0.5mu . \mkern1mu 
\VEC{c}^{\mkern0.5mu n}) = 0$ for $k \ge 1$ and $\vert \mkern1mu 
\textrm{t}_1 \vert + \vert \mkern1mu \textrm{t}_2 \vert + \DOTS + \vert 
\mkern1mu \textrm{t}_k \vert + n < p$.

\smallskip

For any rooted tree $\textrm{t}$, let a row vector $\VEC{D}(\textrm{t})$ 
be defined as $D_{\mkern-0.5mu j \mkern0.5mu}(\textrm{t}) = \bigl( 
\VEC{b} \mkern1mu . \mkern0.5mu \VEC{\Phi}^{\textrm{T}}(\textrm{t}) 
\bigr) \mkern1.5mu \VEC{a}_{* \mkern-0.5mu j} - b_j (1 - 
\smash{c_{\mkern-1.25mu j}^{\ORDT}} ) / \textrm{t}!$ or 
$\VEC{D}(\textrm{t}) = \bigl( \VEC{b} \mkern1mu . \mkern0.5mu 
\VEC{\Phi}{}^{\textrm{T}} \mkern-1mu (\textrm{t}) \bigr) \mkern0.5mu 
\MAT{A} - \bigl( \VEC{b} \mkern0.5mu . \mkern0.5mu (\VEC{1} - 
\VEC{c}^{\ORDT}){}^{\textrm{T}} \bigr) / \textrm{t}!$. For a rooted tree 
$\textrm{t} = \TBPN$ of height at most $1$ the row vector is $\VEC{D} 
\bigl( \TBPN \bigr) = \VEC{d}_n = (\VEC{b} \mkern1mu . \mkern1mu 
\VEC{c}^{\mkern0.5mu n \textrm{T}}) \MAT{A} - \frac{1}{n + 1} \VEC{b} 
\mkern1mu . (\VEC{1} - \VEC{c}^{\mkern0.5mu n + 1}){}^{\textrm{T}}$.

The properties $B(n)$, $C(n)$, and $D(n)$ (for their definition and also 
for simplifying assumptions see, \EG, \cite[p.~52]{But64a}, 
\cite[s.~321]{But16}, \cite[pp.~175, 182, and 208]{HNW93}) can be 
formulated in terms of vectors $\VEC{q}_k$ and $\VEC{d}_k$:

\vspace{-17.5pt}

 \begin{align*}
    B(n): \qquad & \BCN{k} = {\TS\frac{1}{k + 1}} \mbox{ for all }0 \le k 
< n \\[-1pt]
    C(n): \qquad &\VEC{q}_0 = \VEC{q}_1 = \DOTS = \VEC{q}_{n - 1} = 
\VEC{0} \\
    D(n): \qquad &\VEC{d}_0 = \VEC{d}_1 = \DOTS = \VEC{d}_{n - 1} = 
\VEC{0}
 \end{align*}

\vspace{-4.5pt}

Given rooted trees $\textrm{t}_1$, $\textrm{t}_2$, \DOTS, $\TTN$, a new 
tree $\textrm{t}_1 \mkern-1mu \cdot \mkern1mu \textrm{t}_2 \cdot \DOTS 
\cdot \TTN$ is obtained by merging their $n$ roots into one vertex, the 
latter becomes a new root. One has $\VEC{\Phi}(\textrm{t}_1 \mkern-1mu 
\cdot \mkern1mu \textrm{t}_2 \cdot \DOTS \cdot \TTN) = 
\VEC{\Phi}(\textrm{t}_1) \mkern1mu . \mkern1mu \VEC{\Phi}(\textrm{t}_2) 
\mkern1mu . \mkern-1mu \cdot \mkern-2mu \cdot \mkern-2mu \cdot \mkern-1mu 
. \mkern1mu \VEC{\Phi}(\TTN)$. Let $\mathcal{D} : \textrm{T}^2 \to 
\mathcal{T}$ be a mapping defined as $\mathcal{D}(\textrm{t} \SEP 
\textrm{t}\mkern0.75mu{}') = \textrm{t} * \TPRIME + 
{\TS\frac{1}{\textrm{t}!}} [ \mkern1.75mu \bullet^{\ORDT} ] \cdot \TPRIME 
- {\TS\frac1{\textrm{t}!}} \mkern1mu \TPRIME$ for any rooted trees 
$\textrm{t}$ and $\TPRIME$, where $\textrm{t}_1 * \textrm{t}_2 = 
\textrm{t}_1 \cdot [ \mkern1.75mu \textrm{t}_2 \mkern0.5mu ]$ is the 
beta-product of trees \citep{But72}, \citep[p.~45]{But21}. The mapping\\ 
$\mathcal{D}$ could be viewed as a collection of mappings\\ 
$\mathcal{D}(\textrm{t} \SEP \cdot) : \textrm{T} \to \mathcal{T}$, with 
$\TPRIME \mapsto \mathcal{D}(\textrm{t} \SEP \TPRIME)$, indexed by\\ a 
tree $\textrm{t}$, which is a certain way to formalize\\ the concept of 
stumps \citep[\\{}s.~2.7]{But21}. The mapping $\mathcal{D}$ in its 
second\\ agrument could be extended to linear\\ combinations of trees by 
linearity:\\ $\mathcal{D} \bigl( \textrm{t} \SEP \sum_{\TPRIME \in 
\textrm{T}} a(\TPRIME) \mkern1mu \TPRIME \bigr) = \sum_{\TPRIME \in 
\textrm{T}} a(\TPRIME) \mathcal{D}(\textrm{t} \SEP \TPRIME)$.\\ One has 
$\VEC{D}(\textrm{t}) \mkern0.5mu \VEC{\Phi}({\footnotesize\textsf{T}}) = 
\BPHI \bigl( \mathcal{D}(\textrm{t} \SEP {\footnotesize\textsf{T}}) 
\bigr)$ for any $\textrm{t} \in \textrm{T}$\\ and 
${\footnotesize\textsf{T}} \in \mathcal{T}$. The diagram on the right is 
commutative.

\vspace{-122pt}\hfill\begin{picture}(137,111)(0,0)
  \put(34,4){\makebox(0,0)[x]{$\R^s$}}
  \put(94,4){\makebox(0,0)[x]{$\mathcal{T}$}}
  \put(4,56){\makebox(0,0)[x]{$\R$}}
  \put(64,56){\makebox(0,0)[x]{$\R^s$}}
  \put(124,56){\makebox(0,0)[x]{$\mathcal{T}$}}
  \put(94,108){\makebox(0,0)[x]{$\mathcal{T}$}}
  \put(37,61){\makebox(0,0)[x]{\small$\VEC{b}$}}
  \put(96,61){\makebox(0,0)[x]{\small$\VEC{\Phi}$}}
  \put(66,9){\makebox(0,0)[x]{\small$\VEC{\Phi}$}}
  \put(124.5,25){\makebox(0,0)[x]{\small$\mathcal{D}(\textrm{t} \SEP \cdot )$}}
  \put(8,26){\makebox(0,0)[x]{\small$\VEC{D}(\textrm{t})$}}
  \put(112.75,85.5){\makebox(0,0)[x]{\small$\mathcal{Q}$}}
  \put(74,86){\makebox(0,0)[x]{\small$\VEC{Q}$}}
  \put(94,56){\makebox(0,0)[x]{$\xleftarrow{\hspace*{38pt}}$}}
  \put(34,56){\makebox(0,0)[x]{$\xleftarrow{\hspace*{38pt}}$}}
  \put(64,4){\makebox(0,0)[x]{$\xleftarrow{\hspace*{38pt}}$}}
  \put(109,82){\makebox(0,0)[x]{\rotatebox{-60}{$\xrightarrow{\hspace*{38pt}}$}}}
  \put(79,82){\makebox(0,0)[x]{\rotatebox{60}{$\xleftarrow{\hspace*{38pt}}$}}}
  \put(109,30){\makebox(0,0)[x]{\rotatebox{60}{$\xrightarrow{\hspace*{38pt}}$}}}
  \put(19,30){\makebox(0,0)[x]{\rotatebox{-60}{$\xleftarrow{\hspace*{38pt}}$}}}
\end{picture}

\vspace{9.25pt}

The order conditions $\BPHI(\textrm{t}) = 1 / \textrm{t}!$ whenever 
$\ORDT \le p$ could be rewritten as quadrature conditions 
$\BCN{\mkern0.5mu n} = \frac{1}{n + 1}$ for $0 \le n < p$ and order 
conditions of \emph{D-type} $\VEC{D}(\textrm{t}) \mkern0.5mu 
\VEC{\Phi}(\TPRIME) = \nobreak 0$ for all rooted trees $\textrm{t}$ and 
$\TPRIME$ such that $\ORDT + \vert \mkern1mu \TPRIME \vert \le p$. In the 
order conditions of Q-type there could be several $\VEC{Q}(\textrm{t})$ 
vectors, while a condition of D-type has only one $\VEC{D}(\textrm{t})$, 
as a rooted tree can have many branches but only one root.

\section{Dual of a method} \label{dual_method}

\hskip\parindent{}{\bf Definition~\thesection{}.} Consider an $s$-stage 
Runge--Kutta method $\mathscr{M} = (\MAT{A} \SEP \VEC{b} \SEP \VEC{c})$ 
with all the weights being non-zero, \IE, $\BJ \ne 0$ for all $1 \le j 
\le s$, that also satisfies $D(1)$, \IE, $\VEC{b} \mkern1mu \MAT{A} = 
\VEC{b} \mkern1mu . (\VEC{1} - \VEC{c}){}^{\textrm{T}}$. A method 
\emph{dual} to $\mathscr{M}$ is the $s$-stage method $\mathscr{M}^* = 
(\MAT{A}^* \SEP \VEC{b}^* \SEP \VEC{c}^*)$, where $c^*_i = 1 - c_{s + 1 - 
i}$, $b^*_{\mkern-1.25mu j} = b_{s + 1 - j}$, and $a^*_{\mkern0.85mu i 
\mkern-0.85mu j} = b_{s + 1 - j} \mkern1.5mu a_{s + 1 - j \SEP s + 1 - i} 
/ b_{s + 1 - i}$ for all $1 \le i, j \le s$. A method $\mathscr{M}$ is 
called \emph{self-dual} if $\mathscr{M}^* = \mathscr{M}$.

\smallskip

{\bf Statement \thesection{}.1.} Any method $\mathscr{M}$ equals to its 
double dual, \IE, $(\mathscr{M}^*)^* = \mathscr{M}$.

{\bf Statement \thesection{}.2.} If a method is explicit, its dual is 
also explicit.

Examples of self-dual methods are: Kutta's \THIRD{} order method 
\citep[p.~440]{Kut01}, the classic Runge--Kutta method 
\citep[p.~443]{Kut01} and $3/8$ rule \citep[p.~441]{Kut01}, general case 
of $4$-stage methods of order $4$ with symmetrically places nodes 
\citep[eq.~(5.4d)]{But21}, special case of $4$-stage methods of order $4$ 
\citep[p.~442, eq.~(V)]{Kut01} (see also \citep[eq.~(5.4f)]{But21}), 
Gauss--Legendre \citep[]{HaHo55}, \citep[p.~56]{But64a}, Butcher's 
Lobatto \cite[tab.~3]{But64c}, and Lobatto IIIC \citep{Chi71} methods.

\smallskip

A map from a method to its dual is an involution on Runge--Kutta methods 
that exchanges Q- and D-type conditions:

{\bf Theorem \thesection{}.} Consider a method $\mathscr{M}$ that 
satisfies $B(l)$, $C(m)$ and $D(n)$. Then its dual $\mathscr{M}^*$ 
satisfies $B(l)$, $C(n)$ and $D(m)$.

\emph{Proof:} Within the proof $I = s + 1 - i$ and $J = s + 1 - 
\mkern-1mu j$. For any $0 \le k < l$

\vspace{-16pt}

 \begin{align*}
  \VEC{b}^* (\VEC{c}^*)^k &= \sum_{j = 1}^s b_j^* \mkern1mu 
c_{\mkern-1.25mu j}^{* k} = \sum_{j = 1}^s b_J (1 - c_J)^k = \sum_{J = 
1}^s b_J \sum_{k' = 0}^k \mbox{\scalebox{1.2}{$\Bigl($}} \mkern-3mu 
\begin{array}{c} k \\ k' \end{array} \mkern-3mu 
\mbox{\scalebox{1.2}{$\Bigr)$}} (-c_J)^{k'} \\
  &= \sum_{k' = 0}^k (-1)^{k'} \mbox{\scalebox{1.2}{$\Bigl($}} \mkern-3mu 
\begin{array}{c} k \\ k' \end{array} \mkern-3mu 
\mbox{\scalebox{1.2}{$\Bigr)$}} \sum_{J = 1}^s b_J \mkern1mu c_J^{k'} = 
\sum_{k' = 0}^k (-1)^{k'} \mbox{\scalebox{1.2}{$\Bigl($}} \mkern-3mu 
\begin{array}{c} k \\ k' \end{array} \mkern-3mu 
\mbox{\scalebox{1.2}{$\Bigr)$}} \frac{1}{k' + 1} = \frac{1}{k + 1}
 \end{align*}

\vspace{-4pt}

\NI{}Thus $\mathscr{M}^*$ satisfies $B(l)$. For any $0 \le k < n$

\vspace{-16pt}

 \begin{align*}
   \bigl( \MAT{A}^* (\VEC{c}^*)^k \bigr)_i &= \sum_{j = 1}^s 
a_{\mkern1.15mu i \mkern-0.85mu j}^* \mkern1mu c_{\mkern-1.25mu j}^{* k} 
= \sum_{J = 1}^s \frac{b_J \mkern0.5mu a_{JI}}{b_I} (1 - c_J)^k = 
\frac{1}{b_I} \sum_{J = 1}^s b_J \mkern0.5mu a_{JI} \sum_{k' = 0}^k 
\mbox{\scalebox{1.2}{$\Bigl($}} \mkern-3mu \begin{array}{c} k \\ k' 
\end{array} \mkern-3mu \mbox{\scalebox{1.2}{$\Bigr)$}} (-c_J)^{k'} \\
   {} &= \sum_{k' = 0}^k (-1)^{k'} \mbox{\scalebox{1.2}{$\Bigl($}} 
\mkern-3mu \begin{array}{c} k \\ k' \end{array} \mkern-3mu 
\mbox{\scalebox{1.2}{$\Bigr)$}} \frac{1}{b_I} \sum_{J = 1}^s b_J 
\mkern1mu c_J^{k'} a_{JI} = \sum_{k' = 0}^k (-1)^{k'} 
\mbox{\scalebox{1.2}{$\Bigl($}} \mkern-3mu \begin{array}{c} k \\ k' 
\end{array} \mkern-3mu \mbox{\scalebox{1.2}{$\Bigr)$}} \frac{1 - c_I^{k' 
+ 1}}{k' + 1} \\
   {} &= (1 - c_I)^{k + 1} / (k + 1) = c_{\mkern0.5mu i}^{* \SEP k + 1} / 
(k + 1)
 \end{align*}

\vspace{-4pt}

\NI{}Thus $\mathscr{M}^*$ satisfies $C(n)$. For any $0 \le k < m$

\vspace{-16pt}

 \begin{align*}
    \bigr( (\VEC{b}^* \mkern-1mu . \mkern0.5mu (\VEC{c}^{* \textrm{T}}) 
{}^k) \MAT{A}^* \bigr)_{\mkern-1.25mu j} &= \sum_{i = 1}^s b_i^* c_i^{*k} 
a_{ij}^* = \sum_{I = 1}^s b_I (1 - c_I)^k \frac{b_J \mkern0.5mu 
a_{JI}}{b_I} = b_J \sum_{I = 1}^s a_{JI} \sum_{k' = 0}^k 
\mbox{\scalebox{1.2}{$\Bigl($}} \mkern-3mu \begin{array}{c} k \\ k' 
\end{array} \mkern-3mu \mbox{\scalebox{1.2}{$\Bigr)$}} (-c_I)^{k'} \\
    {} &= b_J \sum_{k' = 0}^k (-1)^{k'} \mbox{\scalebox{1.2}{$\Bigl($}} 
\mkern-3mu \begin{array}{c} k \\ k' \end{array} \mkern-3mu 
\mbox{\scalebox{1.2}{$\Bigr)$}} \sum_{I = 1}^s a_{JI} c_I^{k'} = b_J 
\sum_{k' = 0}^k (-1)^{k'} \mbox{\scalebox{1.2}{$\Bigl($}} \mkern-3mu 
\begin{array}{c} k \\ k' \end{array} \mkern-3mu 
\mbox{\scalebox{1.2}{$\Bigr)$}} \frac{c_J^{k' + 1}}{k' + 1} \\
    {} &= b_J \bigl( 1 - (1 - c_J)^{k + 1} \bigr) / (k + 1) = b_j^* 
\bigl( 1 - c_{\mkern-1.25mu j}^{* \SEP k + 1} \bigr) / (k + 1)
 \end{align*}

\vspace{-4pt}

\NI{}Thus $\mathscr{M}^*$ satisfies $D(m)$. \QED

\section{SOL and \BCDQ{} heuristics} \label{heuristics}

\hskip\parindent{}The following is virtually the definition 
\citep[def.~1]{Ver14} of ``\emph{stage order}'':

{\bf Definition~\thesection{}.1.} A stage $i$ is said to be of 
\emph{strong stage order} at least $p$ if $q_{n, i} = \VEC{a}_{i*} 
\mkern0.5mu \VEC{c}^n - \frac{1}{n + 1} c_{\mkern0.5mu i}^{n + 1} = 0$ 
for all $0 \le n < p$, and whenever $\AIJ \ne 0$, then the stage 
$\mkern-0.5mu j$ is of strong stage order at least $p - 1$.

The next definition is weaker but provides more flexibility:

{\bf Definition~\thesection{}.2.} A stage $i$ is said to of \emph{stage 
order} at least $p$ if $Q_i(\textrm{t}) = \VEC{a}_{*i} \mkern1mu 
\VEC{\Phi}(\textrm{t}) - \frac{1}{\textrm{t}!} c_{\mkern0.5mu i}^{\ORDT}
 = 0$ for all rooted trees $\textrm{t}$ with $\ORDT \le p$.

{\bf Statement~\thesection{}.1.} If a stage is of strong stage order 
$p$, then it is of stage order $p$ too.

The advantage of the Definition~\thesection{}.1 is that one can ensure 
that a stage is of a certain stage order, caring only about vectors 
$\VEC{q}_0$, $\VEC{q}_1$, $\VEC{q}_2$, \DOTS{}, instead of 
$\VEC{Q}(\textrm{t})$ for arbitrary rooted trees $\textrm{t}$.

The Stage Order Layers (SOL) heuristic to construct an explicit 
Runge--Kutta method of high order $p$ consists in making stages with 
non-zero weights to be of sufficiently high stage order $q < p$, so that 
many order conditions of Q-type $\VEC{b} \mkern1mu (\VEC{Q}( \mkern0.5mu 
\textrm{t}_1 \mkern-0.5mu ) \mkern0.5mu . \mkern1mu \VEC{Q}( \mkern0.5mu 
\textrm{t}_2 \mkern-0.5mu ) \mkern0.5mu . 
\mkern-2mu\cdot\mkern-2.5mu\cdot\mkern-2.5mu\cdot\mkern-2mu . \mkern1mu 
\VEC{Q}( \mkern0.5mu \textrm{t}_k \mkern-0.5mu ) \mkern0.5mu . \mkern1mu 
\VEC{c}^{\mkern0.5mu n}) = 0$ with $\min \bigl( \vert \textrm{t}_1 \vert 
\SEP \vert \textrm{t}_2 \vert \SEP \DOTS \SEP \vert \textrm{t}_k \vert 
\bigr) \le q$ are automatically satisfied. If $q \ge \frac12 p - 1$, then 
(besides the quadrature conditions) only conditions with $k = 1$, \IE, 
$\VEC{b} \mkern1mu (\VEC{Q}(\textrm{t}) . \mkern1mu \VEC{c}^n) = 0$ with 
$\vert \textrm{t} \vert + n < p$, still need to be dealt with. Sequential 
stage order layers, \IE, subsets of stages with the same stage order, 
greatly increase the redundancy of order conditions of Q-type outside 
(\IE, later stages) of each layer. A prime example of such an approach is 
the \citep{Cur75} method, see Figure~\ref{NCstructure}.

\smallskip

The following definition may be seen as a dual version of the 
Definition~\thesection{}.2:

{\bf Definition~\thesection{}.3.} A stage $j$ is said to of \emph{weak 
stage co-order} at least $p$ if $D_{\mkern-0.5mu j 
\mkern0.5mu}(\textrm{t}) = \bigl( \VEC{b} \mkern1mu . \mkern0.5mu 
\VEC{\Phi}^{\textrm{T}}(\textrm{t}) \bigr) \mkern1.5mu \VEC{a}_{* 
\mkern-0.5mu j} - b_j (1 - \smash{c_{\mkern-1.25mu j}^{\ORDT}} ) / 
\textrm{t}! = 0$ for all rooted trees $\textrm{t}$ with $\vert \textrm{t} 
\vert \le p$.

{\bf Statement~\thesection{}.2.} If a method satisfies $D(1)$, \IE, 
$\VEC{d}_0 = \VEC{b} \mkern1mu \MAT{A} - \VEC{b} \mkern1mu . (\VEC{1} - 
\VEC{c}){}^{\textrm{T}} = \VEC{0}$, then all stages are of weak stage 
co-order at least $1$.

For explicit methods the Definition~\thesection{}.3 is not practically 
useful because it is hardly possible to make even some stages with 
non-zero weights to be of weak stage co-order higher than $1$: Consider 
an explicit $s$-stage method with $\VEC{d}_0 = \VEC{0}$. As $b_s \ne 0$ 
and $d_{\mkern0.5mu 0 \SEP s} = -b_s (1 - c_s) = 0$, one must have $c_s = 
1$.\footnote{~If $c_s = 1$, then $D_s(\textrm{t}) = 0$ for all rooted 
trees $\textrm{t}$, \IE, the stage $s$ has infinite weak stage co-order. 
This is dual to the statement that the stage $1$ has infinite stage order 
as $c_1 = 0$ and $Q_1(\textrm{t}) = 0$ for all $\textrm{t}$.} Then 
$d_{\mkern0.5mu 0 \SEP s - 1} = b_s \mkern1mu a_{s \SEP s - 1} - b_{s - 
1} (1 - c_{s - 1}) = 0$ implies $a_{s \SEP s - 1} = b_{s - 1} (1 - c_{s - 
1}) / b_s$ and $d_{\mkern0.5mu 1 \SEP s - 1} = b_s \mkern1mu a_{s \SEP s 
- 1} - b_{s - 1} (1 - c_{s - 1}^2) / 2 = b_{s - 1} (1 - c_{s - 1})^2 / 
2$. In order to have the stage $(s - 1)$ to be of weak stage co-order at 
least $2$, one must have $a_{s \SEP s - 1} = 0$.\footnote{~This is dual 
to the statement that $q_{1 \SEP 2} = -c_2^2 / 2 = 0$ only if $c_2 = 
a_{21} = 0$.}

\smallskip

{\bf Definition~\thesection{}.4.} Let $\mathscr{S} = \{ 1 \SEP 2 \SEP 
\DOTS \SEP s \}$ be the set of all stages. A \emph{node cluster} is a 
triple $\mathscr{C} = (S \SEP Q \SEP D)$, where $S$ is a non-empty subset 
of $\mathscr{S}$ such that the nodes corresponding to any two stages $i$, 
$j$ in $S$ are identical: $\CI = \CJ$; and $Q \subseteq \R^{\vert S 
\vert}$ and $D \subseteq (\R^{\vert S \vert})^*$ are subspaces (here 
$\vert S \vert$ stands for the cardinality of $S$, and vectors are 
indexed by the elements of $S$) that satisfy the following orthogonality 
conditions:

\vspace{-16pt}

 \begin{align*}
   Q &= \bigl\{ \mkern1mu \VEC{q} \in \R^{\vert S \vert} \,\big\vert\; 
{\TS\sum_{\mkern1mu i \in S}} \mkern1.5mu \BI \mkern1mu q_i = 
{\TS\sum_{\mkern1mu i \in S}} \mkern1.5mu d_i \mkern1mu q_i = 0\mbox{ 
for all }\VEC{d} \in D \mkern1mu \bigr\} \\[1pt]
   D &= \bigl\{ \mkern1mu \VEC{d} \in (\R^{\vert S \vert})^* 
\,\big\vert\; {\TS\sum_{\mkern1mu i \in S}} \mkern1.5mu d_i = 
{\TS\sum_{\mkern1mu i \in S}} \mkern1.5mu d_i \mkern1mu q_i = 0\mbox{ 
for all }\VEC{q} \in Q \mkern1mu \bigr\}
 \end{align*}\IE, $Q$ and $D$ are the orthogonal complements of $D + 
\textrm{span}(\VEC{b} \vert_S)$ and $Q + \textrm{span}(\VEC{1} \vert_S)$, 
respectively. If $\sum_{\mkern1.5mu i \in S} \BI \ne 0$, then 
$\mathscr{C}$ is said to be a \emph{quadrature cluster}.

{\bf Theorem~\thesection{}.} Let $\mathscr{C} = (S \SEP Q \SEP D)$ be a 
node cluster. If $\mathscr{C}$ is a quadrature cluster, then $\VEC{1} 
\vert_S \notin Q$, $\VEC{b} \vert_S \notin D$, and $\dim Q + \dim D = 
\vert S \vert - 1$. If $\mathscr{C}$ is a non-quadrature cluster, \IE, 
$\sum_{\mkern1.5mu i \in S} \BI = 0$, then $\VEC{1} \vert_S \in Q$, 
$\VEC{b} \vert_S \in D$, and $\dim Q + \dim D = \vert S \vert$.

\emph{Proof:} As $Q$ is the orthogonal complement of $D + 
\textrm{span}(\VEC{b} \vert_S)$, its dimension is 
 $\dim Q = 
\vert S \vert - \dim \bigl( D + \textrm{span}(\VEC{b} \vert_S) \bigr)$. 
If $\VEC{b} \vert_S \VEC{1} \vert_S = \sum_{i \in S} \BI \ne 0$, then 
$\VEC{b} \vert_S$ is not in $D$, so $\dim \bigl( D + 
\textrm{span}(\VEC{b} \vert_S) \bigr) = \dim D + 1$.  In a quadrature 
cluster $\VEC{1} \vert_S$ is not orthogonal to $\VEC{b} \vert_S$, and 
thus is not in $Q$. In the case of a non-quadrature cluster, $\sum_{i 
\in S} b_i = 0$, the row vector $\VEC{b} \vert_S$ is orthogonal to $Q + 
\textrm{span}(\VEC{1} \vert_S)$, so $\VEC{b} \vert_S \in D$ and $\dim 
\bigl( D + \textrm{span}(\VEC{b} \vert_S) \bigr) = \dim D$. As the 
vector $\VEC{1} \vert_S$ is orthogonal to both $D$ and $\VEC{b} 
\vert_S$, it is in $Q$. \QED

{\bf Definition~\thesection{}.5.} A node cluster $\mathscr{C} = (S \SEP Q 
\SEP D)$ is said to be of {\it cluster order} at least $p$ if 
$\VEC{Q}(\textrm{t})$ restricted to $S$ is in the subspace $Q$ for all 
rooted trees $\textrm{t}$ with $\ORDT \le p$.

\smallskip

Consider the following two filtrations on the algebra of column vectors 
$\R^s$ with the product taken element-wise \citep[pp.~560 and 
561]{Kha09}, \citep[pp.~683 and 684]{Kha13}:

\vspace{-21.5pt}

 \begin{gather*} \hskip66pt \begin{array}{ccccccccccccc}
   \Phi_0 & \subseteq & \Phi_1 & \subseteq & \Phi_2 & \subseteq & \DOTS & 
\subseteq & \Phi_p & \subseteq & \DOTS & \subseteq & \R^s \\
   \mkern-4mu\mbox{\raisebox{-1.75pt}{\rotatebox{90}{$\subset$}}} & &
   \mkern-4mu\mbox{\raisebox{-1.75pt}{\rotatebox{90}{$\subset$}}} & &
   \mkern-2mu\mbox{\raisebox{-1.75pt}{\rotatebox{90}{$\subseteq$}}} & & & &
   \mkern-2mu\mbox{\raisebox{-1.75pt}{\rotatebox{90}{$\subseteq$}}} & & & &
   \mkern-4mu\mbox{\raisebox{-1.75pt}{\rotatebox{90}{$=$}}} \\
   Q_0 & = & Q_1 & \subseteq & Q_2 & \subseteq & \DOTS & 
\subseteq & Q_p & \subseteq & \DOTS & \subseteq & \R^s
 \end{array} \end{gather*}

\vspace{-2.5pt}

\NI{}The subspace $\Phi_p \subseteq \R^s$ is spanned by the vectors 
$\VEC{\Phi}(\textrm{t})$ for all rooted trees $\textrm{t}$ with $\ORDT 
\le p + 1$. \EEG, $\Phi_0 = \textrm{span}(\VEC{1})$, $\Phi_1 = 
\textrm{span}(\VEC{1} \SEP \VEC{c})$, and $\Phi_2 = \textrm{span}(\VEC{1} 
\SEP \VEC{c} \SEP \VEC{c}^2 \SEP \MAT{A} \mkern0.5mu \VEC{c})$. A 
recursive definition of the subspaces $\Phi_p$ is the following: $\Phi_0 
= \textrm{span}(\VEC{1})$ and $\Phi_p$ is generated by subsets $\Phi_{p - 
1}$, $\MAT{A} \mkern0.5mu \Phi_{p - 1}$, and element-wise products of 
subspaces $\Phi_q \mkern1mu . \mkern0.5mu \Phi_{p - q}$, where $0 < q < 
p$. (For sets $X$ and $Y$ and a binary operation $\star$ the set 
operation is defined as $X \star Y = \{ x \star y \,\vert\, x \in X 
\mbox{~and~} y \in Y \}$, also $x \star Y = \{ x \} \star Y$.) Similarly, 
$Q_0 = Q_1 = \{ \VEC{0} \}$, $Q_2 = \textrm{span}(\VEC{q}_1)$, and $Q_p$ 
is generated by $\VEC{q}_{p - 1}$, $\MAT{A} \mkern0.5mu Q_{p - 1}$, and 
$Q_{\mkern0.5mu q} \mkern1mu . \mkern0.5mu \Phi_{p - q}$, where $1 < q < 
p$. For a method of order $p$ the Q-type order conditions could be 
written as $\VEC{b} \mkern1mu Q_{p - 1} = \{ 0 \}$.

\smallskip

Consider the following non-strictly increasing sequence of subspaces 
(with no structure that is related to their element-wise products) of the 
vector space of row vectors $(\R^s)^*$: $D_{\mkern0.5mu 0} = \{ \VEC{0} 
\}$, $D_{\mkern0.5mu 1} = \textrm{span}(\VEC{d}_0)$, and $D_{\mkern-0.5mu 
p}$ is generated by $D_{\mkern-0.5mu p - 1} . \mkern1mu 
\Phi_1^{\textrm{T}}$, $D_{\mkern-0.5mu p - 1} \mkern1mu \MAT{A}$, and 
$\VEC{D}(\textrm{t})$ for all rooted trees $\textrm{t}$ with $\ORDT = p$. 
For example, if $\VEC{d}_0 = \VEC{0}$, then $\SUBD{2} = 
\textrm{span}(\DO)$ and $\SUBD{3} = \textrm{span} \bigl( \DO \SEP \DO . 
\mkern1mu \VEC{c}^{\textrm{T}} \SEP \DOA \SEP \VEC{d}_2 \SEP (\VEC{b} 
\mkern1mu . (\MAT{A} \mkern0.5mu \VEC{c})^{\textrm{T}}) \mkern0.5mu 
\MAT{A} - \VEC{b} \mkern1mu . (\VEC{1} - \VEC{c}^3)^{\textrm{T}} / 6 
\bigr)$. For a method of order $p$ the D-type order conditions could be 
written as $D_{\mkern-0.5mu p - 1} \VEC{1} = \{ 0 \}$ or, equivalently, 
$D_{\mkern0.25mu q - 1} \mkern1mu \Phi_{p - q} = \{ 0 \}$ for all $1 < q 
\le p$.

A practical working definition of co-order is:

{\bf Definition~\thesection{}.6.} A node cluster $\mathscr{C} = (S \SEP Q 
\SEP D)$ is said to be of {\it cluster co-order} at least $p$ if for any 
row vector $\VEC{d}$ in the subspace $D_{\mkern-0.5mu p}$ its restriction 
to $S$ lies in $D$.

The Counterpoised Node Clusters (CNC) heuristic consists in partitioning 
the set of stages $\mathscr{S}$ into node clusters, and making the node 
clusters order \emph{and} co-order sufficiently high. The word 
``counterpoised'' emphasizes that the subspace $D$ is orthogonal to the 
vector $\VEC{1} \vert_S$ for any node cluster $(S \SEP Q \SEP D)$, that 
helps to satisfy order conditions such as $\VEC{d}_n \mkern1mu 
\VEC{c}^{\mkern0.5mu m} = 0$. Notable examples of using repeated nodes to 
form high order{\SLASH}co-order node clusters are 
\citep[tab.~3.3]{Ver69}, \citep[tab.~1]{CoVe72}, and \citep{Hai78}, see 
Figure~\ref{NCstructure}.

\smallskip

Below the type of a quadrature{\SLASH}non-quadrature node cluster $(S 
\SEP Q \SEP D)$ will be written as $\smash{\CIRCLE{}^{\mkern2mu 
d}_{\mkern2mu d^*}}${\,/\,}$\smash{\Circle{}^{\mkern2mu d}_{\mkern2mu 
d^*}}$ to indicate the dimensions $d = \dim Q$ and $d^* = \dim D$, or as 
$\smash{{}^{\mkern2mu \smash{p}}_{\smash{p^*}}{} \mkern-1.5mu 
\CIRCLE{}^{\mkern2mu \smash{d}}_{\mkern2mu 
\smash{d^*}}}${\,/\,}$\smash{{}^{\mkern2mu \smash{p}}_{\smash{p^*}}{} 
\mkern-1.5mu \Circle{}^{\mkern2mu \smash{d}}_{\mkern2mu \smash{d^*}}}$ to 
additionally indicate its cluster order $p$ and cluster co-order $p^*$.

For a node cluster with only one stage, the possible types are 
$\smash{\CIRCLE{}^{\mkern1.5mu 0}_{\mkern1.5mu 0}}$ and 
$\smash{\Circle{}^{\mkern1mu 1}_{\mkern1.5mu 0}}$. As $\dim D = 0$, and 
the dimension of $Q$ is uniquely determined by whether the cluster is a 
quadrature one or not, the dimensions of $Q$ and $D$ are omitted in one 
node clusters in Figure~\ref{NCstructure}.  With two stages, the types 
can be $\smash{\CIRCLE{}^{\mkern1.5mu 0}_{\mkern1mu 1}}$, 
$\smash{\CIRCLE{}^{\mkern1mu 1}_{\mkern1.5mu 0}}$, 
$\smash{\Circle{}^{\mkern1mu 1}_{\mkern1mu 1}}$, and 
$\smash{\Circle{}^{\mkern1.5mu 2}_{\mkern1.5mu 0}}$. For a node cluster 
of type $\smash{\CIRCLE{}^{\mkern1mu 1}_{\mkern1mu 1}}$ there is an 
interesting possibility: $Q = \textrm{span} \bigl( [ \; g_1 ~\,~ g_2 ~\,~ 
g_3 \; ]{}^{\textrm{T}} \bigr)$, $D = \textrm{span} \bigl( [ \; (g_2 - 
g_3) ~~ (g_3 - g_1) ~~ (g_1 - g_2) \; ] \bigr)$, and $\VEC{b} \vert_S \in 
\textrm{span} \bigl( [ \; (g_2 - g_3) / g_1 ~\,~ (g_3 - g_1) / g_2 ~\,~ 
(g_1 - g_2) / g_3 \; ] \bigr)$. Not only $\VEC{b} \vert_S \mkern1.5mu Q = 
\{ 0 \}$, but $\VEC{b} \vert_S \mkern1mu (Q \mkern1.5mu . \mkern1.5mu Q) 
= \{ 0 \}$ too.

 \begin{figure} \centerline{\includegraphics{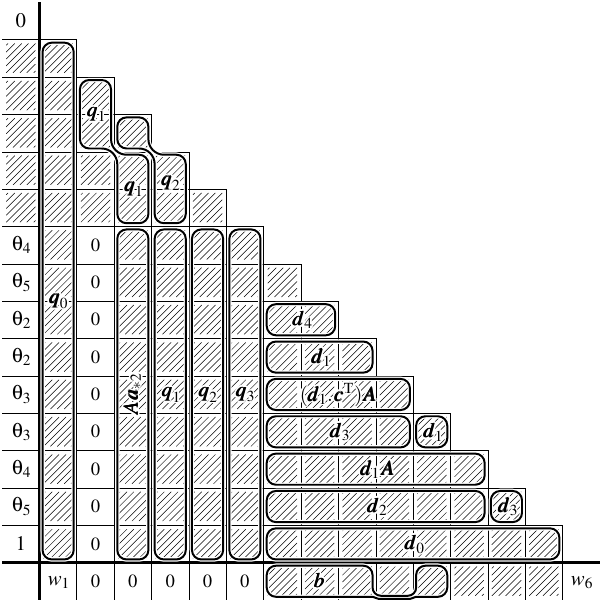}} \vspace{3pt} 
\caption{Butcher tableau \citep[p.~191]{But64b} with all but the final 
steps of the method construction. The shaded cells correspond to the 
nodes, weights, and coefficients that are to be determined. The entries 
marked ``$\VEC{b}$'' are found according to eq.~(\ref{cond_b}); the 
entries marked ``$\VEC{q}_0$'' --- according to eq.~(\ref{cond_q0}); 
``$\VEC{q}_1$'' and ``$\VEC{q}_2$'' in the upper left corner --- 
eq.~(\ref{cond_q1_q2}); parallel blocks with ``$\MAT{A} \mkern1mu 
\VEC{a}_{*2}$'', ``$\VEC{q}_1$'', ``$\VEC{q}_2$'', and ``$\VEC{q}_3$'' 
--- eq.~(\ref{cond_Aa2_q1_q2_q3}); ``$\VEC{d}_0$'' --- 
eq.~(\ref{cond_d0}); ``$\VEC{d}_1$'' --- eq.~(\ref{cond_d1}); 
``$\VEC{d}_2$'' and ``$\DOA$'' --- eq.~(\ref{cond_d2_d1A}); 
``$\VEC{d}_3$'' and ``$(\VEC{d}_1 . \mkern1mu \VEC{c}^{\textrm{T}}) 
\mkern0.5mu \MAT{A}$'' --- eq.~(\ref{cond_d3_cd1A}); and ``$\VEC{d}_4$'' 
--- eq.~(\ref{cond_d4}).} \label{basic_structure} \end{figure}

\section{A family of methods of order 10} \label{family}

\hskip\parindent{}A Runge--Kutta method of order $10$ satisfies order 
conditions $\BPHI(\textrm{t}) = {\TS\frac{1}{\textrm{t}!}}$ for the 
$1205$ rooted trees $\textrm{t}$ such that $\vert \textrm{t} \vert \le 
10$. For a tree $\textrm{t} = [ \mkern1mu \textrm{t}_1 \mkern2mu 
\textrm{t}_2 \mkern2mu \DOTS \mkern2mu \TTN \mkern1mu ]$ the 
corresponding integer partition is $\vert \textrm{t} \vert - 1 = \vert 
\textrm{t}_1 \vert + \vert \textrm{t}_2 \vert + \DOTS + \vert \TTN 
\vert$, see \cite[pp.~50 and 65]{But21}. The partition for $\textrm{t}_1 
\cdot \textrm{t}_2 \cdot \DOTS \cdot \TTN$ is the sum of partitions for 
$\textrm{t}_1$, $\textrm{t}_2$, \DOTS, $\TTN$.

Explicit methods constructed below are based on the $6$-points Lobatto 
quadrature:

\vspace{-16pt}

 \begin{gather*}
   \int\limits_0^1 \mkern-2mu \d\theta \DWS f(\theta) \approx \sum_{k = 
     1}^6 w_k f(\theta_k) , \qquad \alpha, \beta = \sqrt{ 
     \vphantom{\vert^+} \smash{{\TS\frac{1}{21}} (7 \pm 2 \sqrt{7})} } 
     \\
   \theta_1 = 0, \qquad \theta_{2, 5} = {\TS\frac12} (1 \mp \alpha), 
     \qquad \theta_{3, 4} = {\TS\frac12} (1 \mp \beta), \qquad
     \theta_6 = 1 \\
   w_1 = w_6 = {\TS\frac{1}{30}}, \quad w_2 = w_5 = {\TS\frac{1}{60}} 
     (14 - \sqrt{7}), \quad w_3 = w_4 = {\TS\frac{1}{60}} (14 + 
     \sqrt{7})
 \end{gather*}

\vspace{-5pt}

\NI{}The column vector of nodes $\VEC{c}$ is chosen by setting $c_1 = 
\theta_1 = 0$, $c_7 = c_{13} = \theta_4$, $c_8 = c_{14} = \theta_5$, $c_9 
= c_{10} = \theta_2$, $c_{11} = c_{12} = \theta_3$, and $c_{15} = 
\theta_6 = 1$. To satisfy the quadrature order conditions $\BPHI \bigl( 
\TBPN \bigr) = \BCN{\mkern0.5mu n} = \frac{1}{n + 1}$ corresponding to 
partitions $n = \PONE + \PONE + \nobreak \DOTS \nobreak + \nobreak 
\PONE$, where $0 \le n < 10$, the row vector of weights $\VEC{b}$ is 
determined by $b_1 = w_1 = b_{15} = w_6$, $\BJ = 0$ for all $2 \le j \le 
6$, and

\vspace{-15pt}

 \begin{gather}
   b_9 = w_2 - b_{10} , \qquad b_{11} = w_3 - b_{12} , \qquad b_7 = w_4 - 
b_{13} , \qquad b_8 = w_5 - b_{14} \label{cond_b}
 \end{gather}

\vspace{-5pt}

The first column $\VEC{a}_{*1}$ is determined from $\VEC{q}_0 = \MAT{A} 
\VEC{1} - \VEC{c} = \VEC{0}$:

\vspace{-16pt}

 \begin{gather}
   \AIN{1} = c_i - \sum_{j = 2}^{i - 1} \AIJ , \qquad 2 \le i \le 15 
\label{cond_q0}
 \end{gather} 

\vspace{-9pt}

The coefficients $a_{32}$, $a_{42}$, $a_{43}$, $a_{53}$, $a_{54}$, 
$a_{63}$, and $a_{64}$ are found from $\VEC{q}_1 = -\frac{1}{2} c_2^2 
\mkern0.5mu \VEC{e}_2$ and increasing the redundancy of order conditions 
relation $\VEC{q}_2 \in \textrm{span}(\VEC{q}_1 \SEP \MAT{A} \VEC{q}_1)$:

\vspace{-14.5pt}

 \begin{gather}
   \begin{array}{c}
     {\DS a_{32} = \frac{c_3^2}{2 c_2} , \qquad a_{42} = \frac{c_4^2 (3 
c_3 - 2 c_4)}{2 c_2 c_3} , \qquad a_{43} = \frac{c_4^2 (c_4 - 
c_3)}{c_3^2}} \vphantom{\Big\vert_{\Big\vert}} \\
     {\DS a_{53} = \frac{c_5^2 (3 c_4 - 2 c_5) - a_{52} c_2 (6 c_4 - 4 
c_3)}{6 c_3 (c_4 - c_3)} , \qquad a_{54} = \frac{c_5^2 (2 c_5 - 3 c_3) + 
2 a_{52} c_2 c_3}{6 c_4 (c_4 - c_3)}} \vphantom{\Big\vert_{\Big\vert}} \\
     {\DS a_{63} = \frac{c_6^2 (3 c_4 - 2 c_6) - a_{62} c_2 (6 c_4 - 4 
c_3) + 6 a_{65} c_5 (c_5 - c_4)}{6 c_3 (c_4 - c_3)}} \vphantom{\Big\vert_{\Big\vert}} \\
     {\DS a_{64} = \frac{c_6^2 (2 c_6 - 3 c_3) + 2 a_{62} c_2 c_3 - 6 
a_{65} c_5 (c_5 - c_3)}{6 c_4 (c_4 - c_3)}}
   \end{array} \label{cond_q1_q2}
 \end{gather}

\vspace{-5.5pt}

The coefficients $\AIJ$, where $7 \le i \le 15$ and $3 \le j \le 6$, are 
found from $(\MAT{A}^2 \VEC{q}_1)_{\mkern0.5mu i} = -\frac{1}{2} c_2^2 
(\MAT{A} \VEC{a}_{*2})_{\mkern0.5mu i} = 0$ and $q_{1 \SEP i} = q_{2 \SEP 
i} = q_{3 \SEP i} = 0$:

\vspace{-14.5pt}

 \begin{gather}
   \left[ \begin{array}{c} \AIN{3} \\ \AIN{4} \\ \AIN{5} \\ \AIN{6} 
   \end{array} \right] = \left[ \begin{array}{cccc}
     a_{32} & a_{42} & a_{52} & a_{62} \\
     c_3 & c_4 & c_5 & c_6\vphantom{\rule{1pt}{9pt}} \\
     c_3^2 & c_4^2 & c_5^2 & c_6^2\vphantom{\rule{1pt}{11pt}} \\
     ~{c_3^3}~ & ~{c_4^3}~ & ~{c_5^3}~ & ~{c_6^3}~\vphantom{\rule{1pt}{11pt}}
   \end{array} \right]^{-1} \left[ \begin{array}{c}
     0 \\ \frac12 c_i^2 - \sum_{j = 7}^{i - 1} \AIJ \mkern1mu \CJ 
\vphantom{\rule{1pt}{8pt}} \\
     \frac13 c_i^3 - \sum_{j = 7}^{i - 1} \AIJ \mkern1mu c_{\mkern-1.25mu 
j}^2 \vphantom{\rule{1pt}{10pt}} \\
     \frac14 c_i^4 - \sum_{j = 7}^{i - 1} \AIJ \mkern1mu c_{\mkern-1.25mu 
j}^3 \vphantom{\rule{1pt}{10pt}}
   \end{array} \right]
   \label{cond_Aa2_q1_q2_q3}
 \end{gather}

\vspace{-5.5pt}

\NI{}This makes $\VEC{q}_1 = -\frac{1}{2} c_2^2 \mkern0.5mu \VEC{e}_2$, 
$\MAT{A} \VEC{q}_1 = -\frac{1}{2} c_2^2 \mkern0.5mu \VEC{a}_{*2}$, and 
$\VEC{q}_2 = c_2 (c_2 - \frac{2}{3} c_3) \mkern0.5mu \VEC{a}_{*2} - 
\frac{1}{3} c_2^3 \mkern0.5mu \VEC{e}_2$. Also $\VEC{b} \mkern1mu . 
\mkern1mu \VEC{Q}^{\textrm{T}}(\textrm{t}) = \VEC{0}$ for any rooted tree 
$\textrm{t}$ with $\ORDT \le 4$, or any stage $i$ with $b_i \ne 0$ is of 
stage order at least $4$. The order conditions with partitions containing 
parts $\PTWO$, $\PTHREE$, and $\PFOUR$ are satisfied if the ones where 
these parts are decomposed into $\PONE + \PONE$, $\PONE + \PONE + \PONE$, 
and $\PONE + \PONE + \PONE + \PONE$, respectively, are. For any $i$ the 
coefficient $\AIJ$, where $j \le 6$, is now expressed through $c_2$, 
$c_3$, $c_4$, $c_5$, $c_6$, $a_{52}$, $a_{62}$, $a_{65}$, and 
coefficients $a_{\mkern1mu i \mkern1mu k}$ with $k \ge 7$.

The dependence of the coefficients matrix $\MAT{A}$ on the node $c_2$ is 
inconsequential: Due to $\VEC{q}_0 = \VEC{0}$ and $\VEC{q}_1 = 
-\frac{1}{2} c_2^2 \mkern0.5mu \VEC{e}_2$, only the first and second 
columns $\VEC{a}_{*1}$ and $\VEC{a}_{*2}$ depend on $c_2$. The column 
vector $\VEC{a}_{*1} + \VEC{a}_{*2}$ depends on $c_2$ only in its second 
component. The second column $\VEC{a}_{*2}$ is inversely proportional to 
$c_2$.

\smallskip

As $c_{15} = 1$ and $\VEC{a}_{*,15} = \VEC{0}$, one has $d_{n, 15} = 0$ 
for all $n$. The coefficients $a_{15, j}$, where $7 \le j \le 14$, are 
determined from $d_{\mkern0.5mu 0, j} = 0$:

\vspace{-15.5pt}

 \begin{gather}
   a_{15,j} = \frac1{w_6} \biggl( b_j (1 - c_j) - \sum_{i = j + 1}^{14} 
\BI \AIJ \biggr) , \qquad 7 \le j \le 14 \label{cond_d0} \end{gather}

\vspace{-4.5pt}

\NI{}Due to $b_2 = 0$, $\VEC{b} \mkern1mu . \mkern1mu 
\VEC{a}_{*2}^{\textrm{T}} = \VEC{0}$, and $\VEC{b} \mkern1mu . 
\mkern0.5mu (\MAT{A} \VEC{a}_{*2}){}^{\textrm{T}} = \VEC{0}$, one has 
$d_{n,2} = 0$ and $\VEC{d}_n \mkern1mu \VEC{a}_{*2} = 0$ for all $n$. For 
$0 \le m \le 3$, $\VEC{d}_n \mkern1mu \VEC{c}^m = (\VEC{b} \mkern1mu . 
\mkern1mu \VEC{c}^{\mkern0.5mu n\textrm{T}}) \mkern0.5mu \MAT{A} 
\mkern0.5mu \VEC{c}^m - \frac{1}{n + 1} \VEC{b} \mkern0.5mu (\VEC{c}^m - 
\VEC{c}^{m + n + 1}) = (\VEC{b} \mkern0.75mu \VEC{c}^{\mkern0.5mu 
n\textrm{T}}) \VEC{q}_m + \VEC{b} \frac{1}{m + 1} \VEC{c}^{m + n + 1} - 
\frac{1}{(n + 1) (m + 1)} + \frac{1}{(n + 1) (m + n + 2)} = 0 + 
\frac{1}{(m + 1) (m + n + 2)} + \frac{-(m + n + 2) + (m + 1)}{(n + 1) (m 
+ 1) (m + n + 2)} = 0$ for all $n \le 8 - m$. Thus, eq.~(\ref{cond_d0}) 
ensures that $\VEC{d}_0 = \VEC{0}$ and makes the order conditions 
corresponding to partitions $n = n$, where $2 \le n \le 9$, satisfied if 
the ones with the partitions $n - 1 = (n - 1)$ and $n = (n - 1) + \PONE$ 
are.

The order conditions yet to be satisfied could be written as $\VEC{d}_n 
\Phi_m = \{ 0 \}$ conditions of D-type, with $1 \le n \le 4$ and $n + m 
\le 9$. They correspond to partitions $6 = \PFIVE + \PONE$, $7 = \PSIX + 
\PONE$, $8 = \PSEVEN + \PONE$, $9 = \PEIGHT + \PONE$, $7 = \PFIVE + \PONE 
+ \PONE$, $8 = \PSIX + \PONE + \PONE$, $9 = \PSEVEN + \PONE + \PONE$, $8 
= \PFIVE + \PONE + \PONE + \PONE$, $9 = \PSIX + \PONE + \PONE + \PONE$, 
and $9 = \PFIVE + \PONE + \PONE + \PONE + \PONE$.

In order to absorb non-zero values of, \EG, $d_{\mkern0.5mu 1, 14} = w_6 
a_{15, 14} - b_{14} (1 - \theta_5^2) / 2 = b_{14} (1 - \theta_5)^2 / 2$, 
the stages are lumped into node clusters $S_4 = \{ 7 \SEP 13 \}$, $S_5 = 
\{ 8 \SEP 14 \}$, $S_2 = \{ 9 \SEP 10 \}$, and $S_3 = \{ 11 \SEP 12 \}$ 
of type $\smash{{}^{4}_4 \CIRCLE{}^{\mkern1.5mu 0}_{\mkern1mu 1}}$. The 
cluster order $4$ with $Q = \{ \VEC{0} \}$ is already achieved, as all 
the stages from $7$ to $15$ are of stage order $4$. For the cluster 
co-order to be at least $4$, the vectors $\DO$, $\VEC{d}_2$, $\DOA$, 
$\VEC{d}_3$, $(\DO . \mkern1mu \VEC{c}^{\textrm{T}}) \MAT{A}$, $\DO 
\mkern1mu \MAT{A}^2$, and $\VEC{d}_2 \mkern1mu \MAT{A}$ restricted to any 
of these four node clusters should be proportional to the row vector $[ 
\mkern3mu -1 \mkern9mu 1 \mkern3mu ]$. The coefficients $a_{10,7}$, 
$a_{10,8}$, $a_{10,9}$ and $a_{12, 11}$ are found from $d_{\mkern0.5mu 
1,7} + d_{\mkern0.5mu 1,13} = 0$, $d_{\mkern0.5mu 1,8} + d_{\mkern0.5mu 
1,14} = 0$, $d_{\mkern0.5mu 1,9} + d_{\mkern0.5mu 1,10} = 0$, and 
$d_{\mkern0.5mu 1,11} + d_{\mkern0.5mu 1,12} = 0$, respectively:

\vspace{-16pt}

 \begin{align} \begin{array}{rcl} \DS a_{10,7} \mkern-2mu &=& \mkern-2mu 
   {\DS\frac{1}{b_{10} (1 -
      \theta_2)} \biggl( {\TS\frac{1}{2}} w_4 (1 - \theta_4)^2 - 
      \sum_{\substack{i = 8\\i \ne 10}}^{14} \BI (1 - \CI) 
      (a_{\mkern1.15mu i,7} + a_{\mkern1.15mu i,13}) \biggr)} \\
   a_{10,8} \mkern-2mu &=& \mkern-2mu {\DS\frac{1}{b_{10} (1 - \theta_2)} 
     \biggl( {\TS\frac{1}{2}} w_5 (1 - \theta_5)^2 - \sum_{\substack{i = 
     9\\i \ne 10}}^{14} \BI (1 - \CI) (a_{\mkern1.15mu i,8} + 
     a_{\mkern1.15mu i,14}) \biggr)} \\
   a_{10,9} \mkern-2mu &=& \mkern-2mu {\DS\frac{1}{b_{10} (1 - \theta_2)} 
     \biggl( {\TS\frac{1}{2}} w_2 (1 - \theta_2)^2 - \sum_{i = 11}^{14} 
     \BI (1 - \CI) (a_{\mkern1.15mu i,9} + a_{\mkern1.15mu i,10}) 
     \biggr)} \\
   a_{12,11} \mkern-2mu &=& \mkern-2mu {\DS\frac{1}{b_{12} (1 - 
     \theta_3)} \biggl( {\TS\frac{1}{2}} w_3 (1 - \theta_3)^2 - \sum_{i = 
     13}^{14} \BI (1 - \CI) (a_{\mkern1.15mu i,11} + a_{\mkern1.15mu 
     i,12}) \biggr)} \end{array}
   \label{cond_d1}
 \end{align}

\vspace{-2pt}

\NI{}The row vector $\VEC{d}_1$ now has the following structure: 

\vspace{-17pt}

 \begin{align*} \VEC{d}_1 = \bigl[ \mkern3mu 0 \mkern12mu 0 \mkern12mu 0
\mkern12mu 0 \mkern12mu 0 \mkern12mu 0 \mkern9mu -{\mkern-1mu}d_{1,13} 
\mkern9mu -{\mkern-1mu}d_{1,14} \mkern9mu -{\mkern-1mu}d_{1,10} \mkern9mu 
d_{1,10} \mkern9mu -{\mkern-1mu}d_{1,12} \mkern9mu d_{1,12} \mkern9mu 
d_{1,13} \mkern9mu d_{1,14} \mkern12mu 0 \mkern3mu \bigr] \end{align*}

For an order condition $\BPHI(\textrm{t}) = \frac{1}{\textrm{t}!}$, the 
partition of $\textrm{t}$ will be called a $\VEC{b}$-partition. For a 
rooted tree $\textrm{t} = [ \mkern1mu \textrm{t}_1 \mkern2mu \textrm{t}_2 
\mkern2mu \DOTS \mkern2mu \TTM \mkern1mu ]$ with $\ORDT \le 9 - n$, the 
order condition $\VEC{d}_n \VEC{\Phi}(\textrm{t}) = 0$, which is related 
to $\BPHI \bigl( [ \mkern2mu \textrm{t} \mkern4mu \bullet^n \mkern1.5mu ] 
\bigr) = \frac{1}{[ \mkern2mu \textrm{t} \mkern4mu \bullet^n \mkern1.5mu 
]!} = \frac{1}{(\ORDT + n + 1) \mkern2mu \textrm{t}!}$ and corresponds to 
a $\VEC{b}$-partition $\ORDT + n = \ORDT + \PONE + \PONE + \DOTS + 
\PONE$, can be further classified by $\ORDT - 1 = \vert \textrm{t}_1 
\vert + \vert \textrm{t}_2 \vert + \DOTS + \vert \TTM \vert$, which will 
be refferred to as a $\VEC{d}_n$-partition. For any tree $\textrm{t}$ 
with $\ORDT \le 4$, the condition $\DO . \mkern2mu \VEC{Q}{}^{\mkern2mu 
\textrm{T} \mkern-1mu}(\textrm{t}) = \VEC{0}$ holds. Within the order 
conditions $\DO \Phi_8 = \{ 0 \}$, only those with $\VEC{d}_1$-partitions 
having $\PONE$, $\PFIVE$, $\PSIX$, or $\PSEVEN$ as parts need to be 
checked. Due to the structure of the row vector $\DO$, the condition $\DO 
\mkern0.5mu \VEC{c}^n = 0$ with $\DO$-partition $n = \PONE + \PONE + 
\DOTS + \PONE$ (the analogue of a quadrature condition, but for 
$\VEC{d}_1$) is satisfied for all $n$.

The coefficients $a_{14,j}$ and $a_{13,j}$, where $7 \le j \le 12$, are 
found from increasing the redundancy of order conditions relations 

\vspace{-15pt}

 \begin{gather}
   \VEC{d}_2 = \gamma_{20} \mkern1mu \DO + \gamma_{21} (\DO . 
\mkern1mu \VEC{c}^{\textrm{T}}) , \qquad
   \DOA = \gamma_{a0} \mkern1mu \DO + \gamma_{a1} (\DO . 
\mkern1mu \VEC{c}^{\textrm{T}})
 \label{cond_d2_d1A} \end{gather}

\vspace{-5pt}

\NI{}taken at from the $7$th to $12$th components. The four constants 
$\gamma_{20}$, $\gamma_{21}$, $\gamma_{a0}$, and $\gamma_{a1}$ are found 
from these two relations taken at the $13$th and $14$th components. After 
determining $a_{14,j}$ and $a_{13,j}$, where $7 \le j \le 12$, this way 
the first six components of both $\VEC{d}_2$ and of $\DOA$ are equal to 
zero. As $\VEC{d}_0 = 0$, one has $D_{\mkern0.5mu 1} = \{ \VEC{0} \}$, 
$\SUBD{2} = \textrm{span}(\DO)$, and these relations make $\SUBD{3}$ 
being just $2$-dimensional $\textrm{span} \bigl( \DO \SEP \DO . \mkern1mu 
\VEC{c}^{\textrm{T}} \bigr)$. The order conditions corresponding to 
$\VEC{b}$-partitions $n + 2 = (n) + \PONE + \PONE$ are now satisfied if 
the ones with the partitions $n + 1 = (n) + \PONE$ and $n + 2 = (n + 1) + 
\PONE$ are. Also the order conditions with $\DO$-partitions $n = n$, 
where $5 \le n \le 7$, are satisfied if the ones with $n - 1 = (n - 1)$ 
and $n = (n - 1) + \PONE$ are.

The coefficients $a_{11,j}$ and $a_{12,j}$, where $7 \le j \le 10$, and 
the coefficient $a_{14,13}$ are found from increasing the redundancy of 
order conditions relations

\vspace{-15pt}

 \begin{gather}
   \begin{array}{rcl}
     \VEC{d}_3 \mkern-2mu &=& \mkern-2mu {\DS \gamma_{30} \mkern1mu \DO + 
\gamma_{31} (\DO . \mkern1mu \VEC{c}^{\textrm{T}}) + \gamma_{32} (\DO . 
\mkern1mu \VEC{c}^{2\textrm{T}})} \\[4pt]
     (\DO . \mkern1mu \VEC{c}^{\textrm{T}}) \mkern0.5mu \MAT{A} 
\mkern-2mu &=& \mkern-2mu {\DS \gamma_{c0} \mkern1mu \DO + \gamma_{c1} 
(\DO . \mkern1mu \VEC{c}^{\textrm{T}}) + \gamma_{c2} (\DO . \mkern1mu 
\VEC{c}^{2\textrm{T}})}
   \end{array} \label{cond_d3_cd1A}
 \end{gather}

\vspace{-5pt}

\NI{}taken at from the $7$th to $10$th components. The six constants 
$\gamma_{30}$, $\gamma_{31}$, $\gamma_{32}$, $\gamma_{c0}$, 
$\gamma_{c1}$, and $\gamma_{c2}$ are found from these two relations taken 
at the $12$th, $13$th, and $14$th components. The two remaining 
equations, at the $11$th component, are satisfied by tuning the 
coefficient $a_{14,13}$ and having $c_{11} = c_{12}$. (In the $11$-stage 
methods \citep[tab.~3.3, p.~74a]{Ver69} and \citep[tab.~1]{CoVe72} of 
order $8$ the stages $5$, $6$, and $7$ are of stage order $3$; while the 
stages $8$, $9$, $10$, and $11$ are of stage order $4$. In both methods 
$c_7 = c_8$, which is an essential element of the design. Satisfying some 
of the conditions of D-type through the choice of nodes is dual to 
increasing a stage order by setting a node value, see, \EG, 
\citep[eqs.~(6.1), (6.2), (6.3) and (6.4)]{Cur75}.) The first six 
components of both $\VEC{d}_3$ and of $(\DO . \mkern1mu 
\VEC{c}^{\textrm{T}}) \mkern0.5mu \MAT{A}$ are now equal to zero, the 
subspace $\SUBD{4} = \textrm{span} \bigl( \DO \SEP \DO . \mkern1mu 
\VEC{c}^{\textrm{T}} \SEP \DO . \mkern1mu \VEC{c}^{2\textrm{T}} \bigr)$ 
is $3$-dimensional, and the four node clusters based on $S_2$, $S_3$, 
$S_4$, and $S_5$ subsets are of cluster co-order $4$. The order 
conditions corresponding to partitions $8 = \PFIVE + \PONE + \PONE + 
\PONE$ and $9 = \PSIX + \PONE + \PONE + \PONE$ are now satisfied if the 
ones with partitions $6 = \PFIVE + \PONE$, $7 = \PSIX + \PONE$, $8 = 
\PSEVEN + \PONE$, and $9 = \PEIGHT + \PONE$ are. Due to 
eq.~(\ref{cond_d3_cd1A}), verifying the order conditions with 
$\DO$-partitions $6 = \PFIVE + \PONE$ and $7 = \PSIX + \PONE$ is reduced 
to checking the ones with $4 = \PFOUR$, $5 = \PFOUR + \PONE$, $6 = \PFOUR 
+ \PONE + \PONE$, and $5 = \PFIVE$, $6 = \PFIVE + \PONE$, $7 = \PFIVE + 
\PONE + \PONE$, respectively.

The order conditions yet to be satisfied could be written as $(\DO . 
\mkern1mu \VEC{c}^{2\textrm{T}}) \mkern0.5mu \MAT{A} \mkern0.5mu \Phi_4 = 
\{ 0 \}$ and $\VEC{d}_4 \Phi_4 = \{ 0 \}$. They correspond to a 
$\DO$-partition $7 = \PFIVE + \PONE +\PONE$ and $\VEC{d}_4$-partitions $4 
= \PONE + \PONE + \PONE + \PONE$, $4 = \PTWO + \PONE + \PONE$, $4 = \PTWO 
+ \PTWO$, $4 = \PTHREE + \PONE$, $4 = \PFOUR$, respectively, with the 
corresponding $\VEC{b}$-partitions $9 = \PEIGHT + \PONE$ and $9 = \PFIVE 
+ \PONE + \PONE + \PONE + \PONE$. As $\VEC{q}_1 = -\frac12 c_2^2 
\mkern0.5mu \VEC{e}_2$ and $d_{\mkern0.5mu 4,2} = 0$, the 
$\VEC{d}_4$-partitions containing a part $\PTWO$ are reduced to the ones 
where it is decomposed into $\PONE + \PONE$.

The coefficients $a_{97}$ and $a_{98}$ are found from increasing the 
redundancy of order conditions relation

\vspace{-15pt}

 \begin{gather}
   \VEC{d}_4 = \gamma_{40} \mkern1mu \DO + \gamma_{41} (\DO . 
\mkern1mu \VEC{c}^{\textrm{T}}) + \gamma_{42} (\DO . \mkern1mu 
\VEC{c}^{2\textrm{T}}) + \gamma_{43} (\DO . \mkern1mu 
\VEC{c}^{3\textrm{T}}) + \gamma_{4c} (\DO . \mkern1mu 
\VEC{c}^{2\textrm{T}}) \mkern0.5mu \MAT{A}
 \label{cond_d4} \end{gather}

\vspace{-4pt}

\NI{}taken at the $7$th and $8$th components. The five constants 
$\gamma_{40}$, $\gamma_{41}$, $\gamma_{42}$, $\gamma_{43}$, and 
$\gamma_{4c}$ are found from this relation taken at from the $10$th to 
$14$th components. The remaining equation is satisfied by having $c_{9} = 
c_{10}$. The row vectors $\VEC{d}_4$ and $(\VEC{d}_1 . \mkern1mu 
\VEC{c}^{2\textrm{T}}) \mkern0.5mu \MAT{A}$ have their second component 
being equal to zero, but as, \EG, $d_{4,8} + d_{4,14} \ne 0$, their first 
and from the third to sixth components can be non-zero. Nevertheless, the 
relation eq.~(\ref{cond_d4}) is satisfied at all the fifteen components.

\smallskip

It is possible to construct explicit $15$-stage Runge--Kutta methods of 
order $10$ with a different permutation of $6$-points Lobatto quadrature 
nodes. For the design presented here it is necessary that the nodes 
$c_{10}$, $c_{12}$, $c_{13}$, and $c_{14}$ are a permutation of the four 
interior nodes, that $c_9 = c_{10}$ and $c_{11} = c_{12}$, and that the 
stages from $7$ to $14$ use each interior node twice. The nodes in the 
$6$-points Lobatto quadrature are elements of the algebraic extension 
$\FIELD$ of the field of rational numbers $\pmb{\mathbf{Q}}$. Any element 
of $\FIELD$ can be expressed as a linear combination $\xi_1 + \xi_2 
\sqrt{3} + \xi_3 \sqrt{7} + \xi_4 \sqrt{21} + \xi_5 \mkern1mu \alpha + 
\xi_6 \mkern1mu \beta + \xi_7 \sqrt{7} \alpha + \xi_8 \sqrt{7} \beta$ 
with rational weights $\xi_{\mkern1.8mu i \mkern0.5mu}$, $1 \le i \le 8$. 
Such expressions for the fifteen constants $\gamma_{20}$, $\gamma_{21}$, 
$\gamma_{a0}$, $\gamma_{a1}$, $\gamma_{30}$, $\gamma_{31}$, 
$\gamma_{32}$, $\gamma_{c0}$, $\gamma_{c1}$, $\gamma_{c2}$, 
$\gamma_{40}$, $\gamma_{41}$, $\gamma_{42}$, $\gamma_{43}$, and 
$\gamma_{4c}$, which do not depend on $b_{10}$, $b_{12}$, $b_{13}$, and 
$b_{14}$, are given on page~\pageref{constants}.\footnote{~Computations 
were done in interaction with computer algebra system Wolfram 
Mathematica~12.3.0, mainly using commands 
\scalebox{0.85}{\scriptsize\tt\bfseries{}Solve} to symbolically solve 
linear equations, \scalebox{0.85}{\scriptsize\tt\bfseries{}Simplify}, 
\scalebox{0.85}{\scriptsize\tt\bfseries{}Factor}, and 
\scalebox{0.85}{\scriptsize\tt\bfseries{}FindIntegerNullVector}.} The 
list of nine integers $n_1$, $n_2$, \DOTS, $n_9$ corresponds to the 
weights $\xi_{\mkern1.8mu i} = n_{\mkern0.5mu i} / n_9$, $1 \le i \le 8$.

 \begin{figure} 
\centerline{\begin{picture}(327,178)(0,0)
     \linethickness{1.08pt}
     \put(0, 18){\line(1, 0){21}}
     \put(33, 18){\line(1, 0){283}}
     \put(18, 0){\line(0, 1){178}}
     \linethickness{0.36pt}
     \put(36, 0){\line(0, 1){178}}
     \linethickness{0.06pt}
     \put(162, 0){\line(0, 1){178}}
     \put(223, 0){\line(0, 1){158}}
     \put(313, 0){\line(0, 1){138}}
     \put(0, 38){\line(1, 0){21}}
     \put(33, 38){\line(1, 0){284}}
     \put(0, 58){\line(1, 0){21}}
     \put(33, 58){\line(1, 0){284}}
     \put(0, 78){\line(1, 0){21}}
     \put(33, 78){\line(1, 0){284}}
     \put(0, 98){\line(1, 0){21}}
     \put(33, 98){\line(1, 0){284}}
     \put(0, 118){\line(1, 0){21}}
     \put(33, 118){\line(1, 0){284}}
     \put(0, 138){\line(1, 0){21}}
     \put(33, 138){\line(1, 0){281}}
     \put(0, 158){\line(1, 0){21}}
     \put(33, 158){\line(1, 0){190}}
     \put(0, 178){\line(1, 0){21}}
     \put(33, 178){\line(1, 0){129}}
     \put(9.5, 168){\makebox(0, 0)[x]{$\theta_5$}}
     \multiput(9.5, 128)(0, 20){2}{\makebox(0, 0)[x]{$\theta_2$}}
     \multiput(9.5, 88)(0, 20){2}{\makebox(0, 0)[x]{$\theta_3$}}
     \put(9.5, 68){\makebox(0, 0)[x]{$\theta_4$}}
     \put(9.5, 48){\makebox(0, 0)[x]{$\theta_5$}}
     \put(9, 28){\makebox(0, 0)[x]{$1$}}
     \put(28, 128){\makebox(0, 0)[x]{$\cdots$}}
     \put(28, 68){\makebox(0, 0)[x]{$\cdots$}}
     \put(28, 9){\makebox(0, 0)[x]{$\cdots$}}
     \put(99.5, 168){\makebox(0, 0)[x]{$a_{87}$}}
     \put(99.5, 148){\makebox(0, 0)[x]{$\alpha'_{\mkern1mu 24} + u'_2 \mkern0.75mu a_{87}$}}
     \put(99.5, 128){\makebox(0, 0)[x]{$a_{97} + {\TS\frac{w_4{}\vphantom{\vert_\vert}}{b_{10}\vphantom{\vert^\vert}}} (\alpha_{\mkern1mu 24} + u_2 \mkern0.75mu a_{87})$}}
     \put(99.5, 108){\makebox(0, 0)[x]{$\alpha'_{\mkern1mu 34} + u'_3 \mkern0.75mu a_{87}$}}
     \put(99.5, 88){\makebox(0, 0)[x]{$a_{11,7} + {\TS\frac{w_4{}\vphantom{\vert_\vert}}{b_{12}\vphantom{\vert^\vert}}} (\alpha_{\mkern1mu 34} + u_3 \mkern0.75mu a_{87})$}}
     \put(99.5, 68){\makebox(0, 0)[x]{${\TS\frac{w_4{}\vphantom{\vert_\vert}}{b_{13}\vphantom{\vert^\vert}}} (\alpha_{\mkern1mu 44} + u_4 \mkern0.75mu a_{87})$}}
     \put(99.5, 48){\makebox(0, 0)[x]{$a_{87} + {\TS\frac{w_4{}\vphantom{\vert_\vert}}{b_{14}\vphantom{\vert^\vert}}} (\alpha_{\mkern1mu 54} + u_5 \mkern0.75mu a_{87}) - a_{14,13}$}}
     \put(99.5, 28){\makebox(0, 0)[x]{${\TS\frac{w_4{}\vphantom{\vert_\vert}}{w_6{}\vphantom{\vert^\vert}}} (\alpha_{\mkern1mu 64} + u_6 \mkern0.75mu a_{87}) - a_{15,13}$}}
     \put(99.5, 9){\makebox(0, 0)[x]{$w_4 - b_{13}$}}
     \put(192.5, 148){\makebox(0, 0)[x]{$a_{98}$}}
     \put(192.5, 128){\makebox(0, 0)[x]{$a_{98} + {\TS\frac{w_5{}\vphantom{\vert_\vert}}{b_{10}\vphantom{\vert^\vert}}} \mkern0.5mu \alpha_{\mkern1mu 25}$}}
     \put(192.5, 108){\makebox(0, 0)[x]{$a_{11,8}$}}
     \put(192.5, 88){\makebox(0, 0)[x]{$a_{11,8} + {\TS\frac{w_5{}\vphantom{\vert_\vert}}{b_{12}\vphantom{\vert^\vert}}} \mkern0.5mu \alpha_{\mkern1mu 35}$}}
     \put(192.5, 68){\makebox(0, 0)[x]{${\TS\frac{w_5{}\vphantom{\vert_\vert}}{b_{13}\vphantom{\vert^\vert}}} \mkern0.5mu \alpha_{\mkern1mu 45}$}}
     \put(192.5, 48){\makebox(0, 0)[x]{${\TS\frac{w_5{}\vphantom{\vert_\vert}}{b_{14}\vphantom{\vert^\vert}}} \mkern0.5mu \alpha_{\mkern1mu 55}$}}
     \put(192.5, 28){\makebox(0, 0)[x]{${\TS\frac{w_5{}\vphantom{\vert_\vert}}{w_6{}\vphantom{\vert^\vert}}} \mkern0.5mu \alpha_{\mkern1mu 65} - a_{15,14}$}}
     \put(192.5, 9){\makebox(0, 0)[x]{$w_5 - b_{14}$}}
     \put(268.5, 128){\makebox(0, 0)[x]{${\TS\frac{w_2{}\vphantom{\vert_\vert}}{b_{10}\vphantom{\vert^\vert}}} \mkern0.5mu \alpha_{\mkern1mu 22}$}}
     \put(268.5, 108){\makebox(0, 0)[x]{$\alpha'_{\mkern1mu 32} - a_{11,10}$}}
     \put(268.5, 88){\makebox(0, 0)[x]{$\alpha'_{\mkern1mu 32} + {\TS\frac{w_2{}\vphantom{\vert_\vert}}{b_{12}\vphantom{\vert^\vert}}} \mkern0.5mu \alpha_{\mkern1mu 32} - a_{12,10}$}}
     \put(268.5, 68){\makebox(0, 0)[x]{${\TS\frac{w_2{}\vphantom{\vert_\vert}}{b_{13}\vphantom{\vert^\vert}}} \mkern0.5mu \alpha_{\mkern1mu 42} - a_{13,10}$}}
     \put(268.5, 48){\makebox(0, 0)[x]{${\TS\frac{w_2{}\vphantom{\vert_\vert}}{b_{14}\vphantom{\vert^\vert}}} \mkern0.5mu \alpha_{\mkern1mu 52} - a_{14,10}$}}
     \put(268.5, 28){\makebox(0, 0)[x]{${\TS\frac{w_2{}\vphantom{\vert_\vert}}{w_6{}\vphantom{\vert^\vert}}} \mkern0.5mu \alpha_{\mkern1mu 62} - a_{15,10}$}}
     \put(268.5, 9){\makebox(0, 0)[x]{$w_2 - b_{10}$}}
     \put(323, 68){\makebox(0, 0)[x]{$\cdots$}}
     \put(323, 9){\makebox(0, 0)[x]{$\cdots$}}
   \end{picture}}
\vspace{12pt}
\centerline{\begin{picture}(312,118)(0,0)
     \linethickness{1.08pt}
     \put(0, 18){\line(1, 0){21}}
     \put(33, 18){\line(1, 0){268}}
     \put(18, 0){\line(0, 1){118}}
     \linethickness{0.36pt}
     \put(36, 0){\line(0, 1){118}}
     \linethickness{0.06pt}
     \put(110, 0){\line(0, 1){118}}
     \put(173.5, 0){\line(0, 1){98}}
     \put(215, 0){\line(0, 1){78}}
     \put(256.5, 0){\line(0, 1){58}}
     \put(298, 0){\line(0, 1){38}}
     \put(0, 38){\line(1, 0){21}}
     \put(33, 38){\line(1, 0){265}}
     \put(0, 58){\line(1, 0){21}}
     \put(33, 58){\line(1, 0){223.5}}
     \put(0, 78){\line(1, 0){21}}
     \put(33, 78){\line(1, 0){182}}
     \put(0, 98){\line(1, 0){21}}
     \put(33, 98){\line(1, 0){140.5}}
     \put(0, 118){\line(1, 0){21}}
     \put(33, 118){\line(1, 0){77}}
     \multiput(9.5, 88)(0, 20){2}{\makebox(0, 0)[x]{$\theta_3$}}
     \put(9.5, 68){\makebox(0, 0)[x]{$\theta_4$}}
     \put(9.5, 48){\makebox(0, 0)[x]{$\theta_5$}}
     \put(9, 28){\makebox(0, 0)[x]{$1$}}
     \put(28, 68){\makebox(0, 0)[x]{$\cdots$}}
     \put(28, 9){\makebox(0, 0)[x]{$\cdots$}}
     \put(73.5, 108){\makebox(0, 0)[x]{${\TS\frac{b_{10}\vphantom{\vert_\vert}}{w_3{}\vphantom{\vert^\vert}}} \mkern0.25mu \textrm{A}_{11,10}$}}
     \put(73.5, 88){\makebox(0, 0)[x]{$a_{11,10} + {\TS\frac{b_{10}\vphantom{\vert_\vert}}{b_{12}\vphantom{\vert^\vert}}} \mkern0.25mu \textrm{A}_{12,10}$}}
     \put(73.5, 68){\makebox(0, 0)[x]{${\TS\frac{b_{10}\vphantom{\vert_\vert}}{b_{13}\vphantom{\vert^\vert}}} \mkern0.25mu \textrm{A}_{13,10}$}}
     \put(73.5, 48){\makebox(0, 0)[x]{${\TS\frac{b_{10}\vphantom{\vert_\vert}}{b_{14}\vphantom{\vert^\vert}}} \mkern0.25mu \textrm{A}_{14,10}$}}
     \put(73.5, 28){\makebox(0, 0)[x]{${\TS\frac{b_{10}\vphantom{\vert_\vert}}{w_6{}\vphantom{\vert^\vert}}} \mkern0.25mu \textrm{A}_{15,10}$}}
     \put(73.5, 9){\makebox(0, 0)[x]{$b_{10}$}}
     \put(142, 88){\makebox(0, 0)[x]{${\TS\frac{w_3{}\vphantom{\vert_\vert}}{b_{12}\vphantom{\vert^\vert}}} \mkern0.5mu \alpha_{\mkern1mu 33}$}}
     \put(142, 68){\makebox(0, 0)[x]{${\TS\frac{w_3{}\vphantom{\vert_\vert}}{b_{13}\vphantom{\vert^\vert}}} \mkern0.5mu \alpha_{\mkern1mu 43} - a_{13,12}$}}
     \put(142, 48){\makebox(0, 0)[x]{${\TS\frac{w_3{}\vphantom{\vert_\vert}}{b_{14}\vphantom{\vert^\vert}}} \mkern0.5mu \alpha_{\mkern1mu 53} - a_{14,12}$}}
     \put(142, 28){\makebox(0, 0)[x]{${\TS\frac{w_3{}\vphantom{\vert_\vert}}{w_6{}\vphantom{\vert^\vert}}} \mkern0.5mu \alpha_{\mkern1mu 63} - a_{15,12}$}}
     \put(142, 9){\makebox(0, 0)[x]{$w_3 - b_{12}$}}
     \put(194.5, 68){\makebox(0, 0)[x]{${\TS\frac{b_{12}\vphantom{\vert_\vert}}{b_{13}\vphantom{\vert^\vert}}} \mkern0.25mu \textrm{A}_{13,12}$}}
     \put(194.5, 48){\makebox(0, 0)[x]{${\TS\frac{b_{12}\vphantom{\vert_\vert}}{b_{14}\vphantom{\vert^\vert}}} \mkern0.25mu \textrm{A}_{14,12}$}}
     \put(194.5, 28){\makebox(0, 0)[x]{${\TS\frac{b_{12}\vphantom{\vert_\vert}}{w_6{}\vphantom{\vert^\vert}}} \mkern0.25mu \textrm{A}_{15,12}$}}
     \put(194.5, 9){\makebox(0, 0)[x]{$b_{12}$}}
     \put(236, 48){\makebox(0, 0)[x]{${\TS\frac{b_{13}\vphantom{\vert_\vert}}{b_{14}\vphantom{\vert^\vert}}} \mkern0.25mu \textrm{A}_{14,13}$}}
     \put(236., 28){\makebox(0, 0)[x]{${\TS\frac{b_{13}\vphantom{\vert_\vert}}{w_6{}\vphantom{\vert^\vert}}} \mkern0.25mu \textrm{A}_{15,13}$}}
     \put(236, 9){\makebox(0, 0)[x]{$b_{13}$}}
     \put(277.5, 28){\makebox(0, 0)[x]{${\TS\frac{b_{14}\vphantom{\vert_\vert}}{w_6{}\vphantom{\vert^\vert}}} \mkern0.25mu \textrm{A}_{15,14}$}}
     \put(277.5, 9){\makebox(0, 0)[x]{$b_{14}$}}
     \put(308, 9){\makebox(0, 0)[x]{$\cdots$}}
   \end{picture}} \vspace{4pt} \caption{Bottom right corner of the 
Butcher tableau obtained by satisfying the quadrature conditions 
eq.~(\ref{cond_b}) and increading the redundancy in the order conditions 
of D-type relations eqs.~(\ref{cond_d0}), (\ref{cond_d1}), 
(\ref{cond_d2_d1A}), (\ref{cond_d3_cd1A}), and (\ref{cond_d4}). The 
weights $b_{10}$, $b_{12}$, $b_{13}$, and $b_{14}$ are free parameters. 
The coefficient $a_{87}$ is determined later. The exact numerical values 
of constants $\textrm{A}_{15,14}$, $\textrm{A}_{15,13}$, \DOTS, 
$\alpha_{\mkern1mu 32}$, $\alpha'_{\mkern1mu 32}$, \DOTS, $u_2$, and 
$u'_2$ are given on page~\pageref{constants}.} \label{closing} 
\end{figure}

All coefficients $\AIJ$ are now expressed through $c_2$, $c_3$, $c_4$, 
$c_5$, $c_6$, $b_{10}$, $b_{12}$, $b_{13}$, $b_{14}$, $a_{52}$, 
$a_{62}$, $a_{65}$, and $a_{87}$. The structure of the bottom right 
corner of the Butcher tableau, \IE, the coefficients $\AIJ$ for $j \ge 
7$, is shown in Figure~\ref{closing}. The exact expressions for the 
forty two constants $\textrm{A}_{15,14}$, $\textrm{A}_{15,13}$, 
$\textrm{A}_{14,13}$, \DOTS, $\alpha_{\mkern1mu 32}$, 
$\alpha'_{\mkern1mu 32}$, $\alpha_{\mkern1mu 22}$, \DOTS, $u_3$, $u'_3$, 
$u_2$, and $u'_2$ are given on page~\pageref{constants}. When dealing 
with the order conditions of D-type, to at least partially eliminate the 
presence of the weights, it is convenient to use renormalized by weights 
variables $\AIJ = \BJ \mkern1mu \textrm{A}_{\mkern1.25mu i \mkern-0.85mu 
j} / \mkern0.5mu \BI$ and $d_{n \SEP j} = \BJ \mkern1mu 
\Delta_{\mkern1.25mu n \mkern-0.75mu j}$:

\vspace{-15pt}

 \begin{gather*}
   \textrm{A}_{15, j} = 1 - \CJ \,\, - \!\! \sum_{i = j + 1}^{14} \! 
\textrm{A}_{\mkern1.25mu i \mkern-0.85mu j} , \qquad
   \Delta_{\mkern1.25mu n \mkern-0.75mu j} = \DN(\CJ) \,\, - \!\! 
\sum_{i = j + 1}^{14} (1 - \CIN) \textrm{A}_{\mkern1.25mu i \mkern-0.85mu j} 
\end{gather*}

\vspace{-4pt}

\NI{}where $\DN(\theta) = 1 - \theta - {\TS\frac1{n + 1}} (1 - 
\theta^{\mkern1mu n + 1})$, \EG, $\mathcal{D}_{\mkern1.75mu 0}(\theta) = 
0$ and $\mathcal{D}_{\mkern1.75mu 1}(\theta) = {\TS\frac12} (1 - 
\theta)^2$. The coefficients $\textrm{A}_{\mkern1mu i \mkern-0.85mu j}$ 
are of a would-be dual method (as $b_2 = b_3 = b_4 = b_5 = b_6 = 0$, the 
dual method does not exist). \smallskip

The four remaining order conditions to be satisfied are $\VEC{d}_4 
\mkern1mu \VEC{c}^4 = \VEC{d}_4 (\VEC{c} . \mkern1mu \VEC{a}_{*2}) = 
\VEC{d}_4 \mkern1mu \MAT{A} \mkern1mu \VEC{a}_{*2} = \VEC{d}_4 \mkern1mu 
\VEC{q}_3 = 0$. By dimension counting, satisfying them would reduce the 
number of free parameters by four, resulting in a $9$-dimensional family 
of methods of order $10$. Satisfying the remaining conditions with 
maximal possible generality is cumbersome, though. One way to simplify 
the further analysis is to set $c_3 = \theta_3$, $c_4 = \theta_4$, $c_5 = 
\theta_5$, $c_6 = \theta_2$, then it is possible to construct a 
$5$-dimensional family of explicit $15$-stage methods of order $10$, 
parametrized by $c_2$, $b_{10}$, $b_{12}$, $b_{13}$, and $b_{14}$, with 
coefficients in a certain quadratic extension of $\FIELD$.

A more sensible approach would be to increase the number of stage order 
layers in the opening stages. Let $c_3 = \frac{2}{3} c_4$ (which implies 
$a_{42} = 0$) and $a_{52} = a_{62} = 0$. With $\VEC{a}_{*2} = a_{32} 
\mkern1mu \VEC{e}_3$, the coefficients $\AIN{3}$, where $7 \le i \le 15$, 
are all zero. The stages $2$, $3$, from $4$ to $6$, and from $7$ to $15$ 
are of strong stage order $1$, $2$, $3$, and $4$, respectively. With 
$d_{\mkern0.5mu 4,3} = 0$ the condition $\VEC{d}_4 (\VEC{c} . \mkern1mu 
\VEC{a}_{*2}) = 0$ is satisfied. The row vector $\VEC{d}_4$ does not 
depend on $a_{65}$. The coefficient $a_{87}$ is found from the condition 
$\VEC{d}_4 \mkern1mu \VEC{c}^4 = 0$ and is now expressed through $c_4$, 
$c_5$, and $c_6$. The coefficient $a_{65}$ is found from $\VEC{d}_4 
\mkern1mu \MAT{A} \mkern1mu \VEC{a}_{*2} = 0$. The last remaining order 
condition $\VEC{d}_4 \mkern1mu \VEC{q}_3 = 0$ is satisfied by setting the 
node $c_6$:
 \begin{gather*}
   c_6 = \frac {U (c_4 + c_5) + 14 U' c_4 c_5 + U'' c_4 c_5 (c_4 + c_5)} 
{3 U + 14 U' (c_4 + c_5) + 2 U'' (c_4^2 + c_5^2) + 7 c_4 c_5 \bigl( V 
+ 20 V' (c_4 + c_5) + 60 V'' c_4 c_5 \bigr)}
 \end{gather*}The exact values of the constants $U$, $U'$, $U''$, $V$, 
$V'$, and $V''$ are given on page~\pageref{constants}. The result is a 
$7$-dimensional family of explicit $15$-stage Runge--Kutta methods of 
order $10$, parametrized by $c_2$, $c_4$, $c_5$, $b_{10}$, $b_{12}$, 
$b_{13}$, and $b_{14}$. The following choice of parameters gives a 
method with comparatively low magnitude of the coefficients: 
\begin{gather}
    c_2 = {\TS\frac{2}{15}}, \!\quad\! c_4 = {\TS\frac{2}{5}}, \!\quad\! 
c_5 = {\TS\frac{4}{7}}, \!\quad\! b_{10} = {\TS\frac{2}{7}} w_2 , 
\!\quad\! b_{12} = {\TS\frac{2}{9}} w_3, \!\quad\! b_{13} = w_4, 
\!\quad\! b_{14} = w_5 \label{the_method}
 \end{gather} The method is presented on page~\pageref{decimal} in its 
rounded decimal form. The format is $15$ numbers for the nodes $\VEC{c}$, 
$15$ numbers for the weights $\VEC{b}$, and $1 + 2 + 3 + \DOTS + 14 = 
105$ numbers for below the diagonal part of the coefficients matrix 
$\MAT{A}$, row by row. The dataset 
\scalebox{0.95}[1]{\href{https://doi.org/10.5281/zenodo.15429897}{\tt\bfseries{}doi:10.5281/zenodo.15429897}} 
contains the exact values.

Under the condition $c_4 + c_5 = 3 (1 + \beta) / 4$ the expression for 
$c_6$ simplifies to $c_6 = U_1 / \bigl( U_2 - (c_5 - c_4)^2 \bigr)$ for 
some constants $U_1$ and $U_2$.

 \section{Properties of some methods of order 10} \label{comparison}

\hskip\parindent{}The basic properties of some known explicit 
Runge--Kutta methods of order $10$ and of the new method 
eq.~(\ref{the_method}) are compared in Table~\ref{table_methods} and in 
Figures~\ref{NCstructure}, \ref{region}, and \ref{circle}, where the 
methods are named as follows: C10 is \citep{Cur75}; H10 is 
\citep{Hai78}, O10 is \citep{Ono03}; F10 is \citep{Fea07}; and Z10 is 
\citep{Zha24}.

In C10 the stages from $2$ to $11$ are forming five stage order layers 
$\{ 2 \}$, $\{ 3 \}$, $\{ 4 \SEP 5 \}$, $\{ 6 \SEP 7 \}$, and $\{8 \SEP 9 
\SEP 10 \SEP 11 \}$, see Figure~\ref{NCstructure}, top left panel. To 
absorb non-zero values of, \EG, of $d_{1,17}$, virtually out of necessity 
two counterpoised node clusters $\{ 14 \SEP 17 \}$ and $\{ 12 \SEP 16 \}$ 
of type $\smash{{}^{6}_3 \CIRCLE{}^{\mkern1.5mu 0}_{\mkern1mu 1}}$ are 
used. Still, much effort is spent in the opening for the later stages 
(from $12$ to $18$) to be of stage order $6$.

In H10 there are four non-quadrature node clusters $\{ 2 \SEP 16 \}$, $\{ 
3 \SEP 15 \}$, $\{ 6 \SEP 13 \}$, and $\{ 7 \SEP 14 \}$, see 
Figure~\ref{NCstructure}, top right panel. They play a role of nested 
layers of insulation, both from the opening and closing, around the four 
stages $9$, $10$, $11$, and $12$, allowing the latter to have both high 
stage order ($5$) and stage co-order ($4$).

Regions of absolute stability of the six methods are shown in 
Figure~\ref{region}.

The internal structure of the methods through the progression, 
sensitivity to the r.h.s. function, and alignment along the trajectory of 
the intermediate positions $\X12S{i}\mkern0.5mu$, $1 \le i \le s$, is 
shown in Figure~\ref{circle}.

\newpage

\begin{table} \centerline{\small\begin{tabular}{r|c|ccc|cl}
     & $s$ & $10^6 \times T_{11}$ & $10^6 \times T_{12}$ & $10^6 \times 
T_{13}$ & {$\phantom{-}\mbox{max}_{ij} \vert \mkern1mu \AIJ \vert$} & 
{\,\,$\phantom{-}\mbox{min}_j \mkern1.5mu \BJ \STRUT$} \\
     \hline
     C10$\STRUT$ & \,$18$\, & $\phantom{0}3.50...$ & $\phantom{0}8.14...$ & 
$\phantom{0}13.06...$ & $\phantom{-}5.4724...$ & $\phantom{-}0.03333...$ \\
     H10 & $17$ & $\phantom{0}5.27...$ & $17.22...$ & 
$\phantom{0}36.01...$ & $\phantom{-}1.0549...$ & $-0.18$ \\
     O10$\STRUT$ & $17$ & $\phantom{0}1.25...$ & $\phantom{0}3.01...$ & 
$\phantom{00}4.71...$ & $\phantom{-}1.3763...$ & $-0.17892...$ \\
     F10 & $17$ & $21.89...$ & $64.01...$ & $113.71...$ & $\phantom{-}5.7842...$ & 
$-0.05$ \\
     Z10$\STRUT$ & $16$ & $\phantom{0}1.42...$ & $21.70...$ & 
$\phantom{0}37.89...$ & $\phantom{-}4.9406...$ & $-1.19177...$ \\
     eq.~(\ref{the_method})$\vphantom{\vert^\vert}$ & $15$ & 
$\phantom{0}3.49...$ & $\phantom{0}8.48...$ & \phantom{0}$14.07...$ & 
$\phantom{-}2.2415...$ & $\phantom{-}0.03333...$
 \end{tabular}\hskip18pt} \vspace{6pt} \centerline{\small\begin{tabular}{r|l|ll|ll}
     & $\phantom{-0}z_R$ & \hskip15pt$\tilde{x}(\pi/2)$ & 
\hskip9pt$\tilde{y}(\pi/2)$ & \hskip12pt$\tilde{x}(\pi/2)$ & 
\hskip6pt$\tilde{y}(\pi/2)$ \\
     \hline
     C10$\STRUT$ & $-3.8269...$ & $-0.00001559...$ & $1.0000226...$ & 
$\phantom{-}0.000093...$ & $1.000561...$ \\
     H10 & $-2.7046...$ & $-0.00071183...$ & $1.0004307...$ & 
$\phantom{-}0.011791...$ & $1.007904...$ \\
     O10$\STRUT$ & $-3.3815...$ & $-0.00006422...$ & $1.0000264...$ & 
$\phantom{-}0.000151...$ & $1.000116...$ \\
     F10 & $-2.5279...$ & $-0.00091244...$ & $1.0007372...$ & 
$-0.004805...$ & $0.996073...$ \\
     Z10$\STRUT$ & $-4.7240...$ & $-0.00000464...$ & $1.0000090...$ & 
$-0.004199...$ & $0.997594...$ \\
     eq.~(\ref{the_method})$\vphantom{\vert^\vert}$ & $-4.4293...$ & 
$-0.00000074...$ & $1.0000335...$ & $\phantom{-}0.000203...$ & 
$1.000054...$
 \end{tabular}} \vspace{6pt} \caption{A comparison of six explicit 
$s$-stage Runge--Kutta methods of order $10$. Error coefficients are 
defined as $T_p^2 = \sum_{\;\textrm{t},\;\vert \mkern1mu \textrm{t} 
\mkern1mu \vert = p} \bigl( \VEC{b} \mkern2mu \VEC{\Phi}(\textrm{t}) - 1 
/ \textrm{t}! \bigr){}^2 \mkern-1mu / \mkern1mu \sigma^2(\textrm{t}) = (1 
/ p!)^2 \sum_{\;\textrm{t},\;\vert \mkern1mu \textrm{t} \mkern1mu \vert = 
p} \alpha^2(\textrm{t}) \bigl( \textrm{t}! \mkern3mu \VEC{b} \mkern2mu 
\VEC{\Phi}(\textrm{t}) - 1 \bigr){}^2$, where $\sigma(\textrm{t})$ is the 
order of the symmetry group of the tree $\textrm{t}$, and 
$\alpha(\textrm{t})$ is the number of monotonic labelings of $\textrm{t}$ 
(see, \EG, \citep[ss.~304 and 318]{But16}, \citep[pp.~147 and 
158]{HNW93}, \citep[pp.~58 and 60]{But21}, \citep[pp.~57 and 58]{HLW06}). 
The $\mbox{min}_j \mkern1.5mu \BJ$ column shows the minimal value of a 
non-zero weight. The interval of absolute stability $[ \mkern1.5mu z_R 
\SEP 0 \mkern1.5mu ]$ is a connected component of $\bigl\{ \mkern0.5mu z 
\mkern2mu \big\vert \mkern2mu z \in \R \textrm{ and } \vert R(z) \vert 
\le 1 \bigr\}$ that contains zero, here $R(z) = 1 + \sum_{n = 0}^{s - 1} 
\mkern1mu z^{\mkern1mu n + 1} \VEC{b} \mkern2mu \MAT{A}^n \VEC{1}$ is the 
stability function (see, \EG, \citep[s.~238]{But16}, 
\citep[sec.~4.4]{AsPe98}, \citep[s.~5.3]{But21}); see also 
Figure~\ref{region}. The left and right pairs of columns 
$\tilde{x}(\pi/2)$, $\tilde{y}(\pi/2)$ give the result of the application 
of one step $h = \pi / 2$ to systems of differential equations $\d x / 
\DT = -y$, $\d y / \DT = x$ and $\d x / \DT = -y / (x^2 + y^2)$, $\d y / 
\DT = x / (x^2 + y^2)$, respectively, with the initial condition $x(0) = 
1$, $y(0) = 0$; see also Figure~\ref{circle}. For both systems the exact 
solution is $x\BRAT = \cos t$, $y\BRAT = \sin t$. For the left columns 
pair one has $\tilde{x}(\pi/2) + \i \, \tilde{y}(\pi/2) = R(\i \, \pi / 
2)$.} \label{table_methods} \end{table}

\begin{figure} \centerline{\input{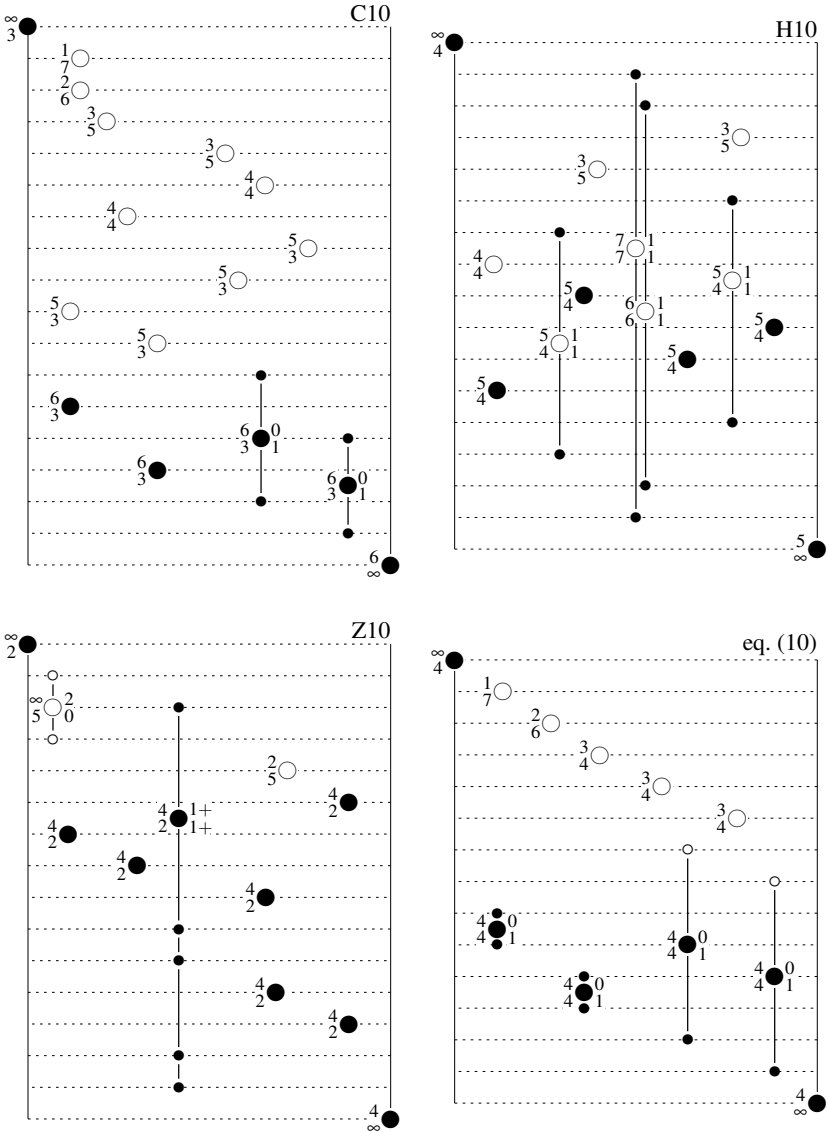}} \vspace{6pt} 
\caption{Node cluster structures of four methods of order $10$: C10, H10, 
Z10, and eq.~(\ref{the_method}). The dashed lines mark the stages, from 
$1$ (upper line) to $s$ (lower line). The horizontal axis represents the 
node positions. Vertical lines connecting small circles correspond to 
node clusters with more than one stage. As the Z10 method was found by 
numerically minimizing $\mathcal{F} = \sum_{\mkern2mu \textrm{t}} \bigl( 
\BPHI(\textrm{t}) - 1 / \textrm{t}! \bigr){}^2$, where the summation goes 
over all rooted trees $\textrm{t}$ with $\vert \textrm{t} \vert \le 10$, 
with the value $\mathcal{F} = 0$ being eventually achieved; the method 
lacks ``explainability'': The reasons behind its structure are not 
transparent or easily understandable. The Z10 method contains an 
idiosyncratic node cluster $(S \SEP Q \SEP D)$ with $5$ stages: $S = \{ 3 
\SEP 10 \SEP 11 \SEP 14 \SEP 15 \}$. On this cluster $\VEC{b} \vert_S 
\mkern2mu \VEC{q}_4 \vert_S = -1.28... \times 10^{-4} \ne 0$ and 
$\VEC{d}_2 \vert_S \mkern2mu \VEC{1} \vert_S = 1.62... \times 10^{-5} \ne 
0$, and its cluster order and co-order are $4$ and $2$, respectively. 
With $\VEC{q}_1 \vert_S = \VEC{0}$, the vectors $\VEC{q}_2 \vert_S$, 
$(\MAT{A} \VEC{q}_1) \vert_S$, $\VEC{q}_3 \vert_S$, and $(\MAT{A} 
\VEC{q}_2) \vert_S$ generate a one-dimensional subspace $Q' \subseteq 
\R^5$. With $\VEC{d}_0 = \VEC{0}$, the subspace $D' = 
\textrm{span}(\VEC{d}_1 \vert_S)$ is also one-dimensional. The subspaces 
$Q \supseteq Q'$ and $D \supseteq D'$, with $\dim Q + \dim D = 4$, can be 
chosen in a variety of ways (that is why the cluster type is shown as 
${}^4_2 \mkern-0.5mu \CIRCLE \mkern0.5mu {}^{1+}_{1+}$).} 
\label{NCstructure} \end{figure}

\begin{figure} \centerline{\includegraphics{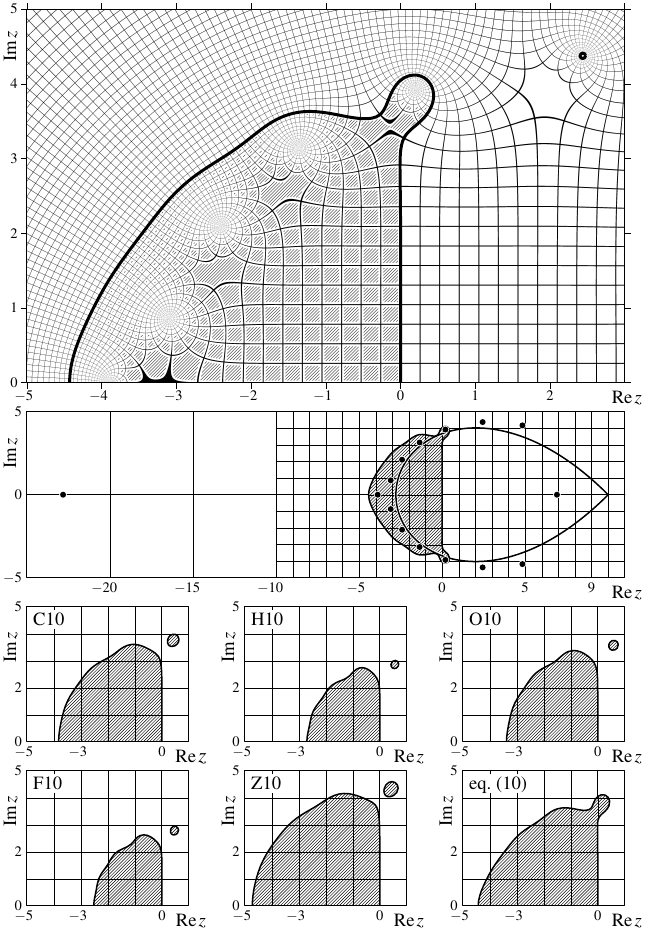}} \vspace{6pt} 
\caption{\looseness=-1{}Regions of absolute stability $\bigl\{ 
\mkern0.5mu z \mkern2mu \big\vert \mkern2mu \vert R(z) \vert \le 1 
\bigr\}$ of C10, H10, O10, F10, Z10, and eq.~(\ref{the_method}) methods. 
The upper and middle panels correspond to the method 
eq.~(\ref{the_method}). In the upper panel the thick solid curve marks 
the boundary of the shaded stability region. The curvilinear grid of 
lines depicts the regions $d \bigl( \textrm{Re} \, L(z) \bigr) \le 1/40$ 
and $d \bigl( \textrm{Im} \, L(z) \bigr) \le 1/40$, where $L(z) = (12 / 
\pi) \log R(z)$ and $d(x) = \min_{n \in \Z} \vert \mkern1.5mu x - n \vert 
= \mbox{$\vert \mkern1.5mu x - \textrm{round}(x) \vert$}$ is the distance 
to the closest integer. The grid becomes more dense on the left due to 
the argument principle and numerous zeros of $R(z)$. As $R(z) \approx 
\exp(z)$ in the vicinity of $z = 0$ the curvilinear grid resembles a 
square one there. In the middle panel the fifteen points correspond to 
zeros of the stability function $R(z)$, while the solid curve is the 
Szeg\H{o} curve $\vert \mkern1mu z \exp(1 - z) \vert = 1$ \citep{Sze24} 
expanded by factor $10$.} \label{region} \end{figure}

\begin{figure} \centerline{\includegraphics{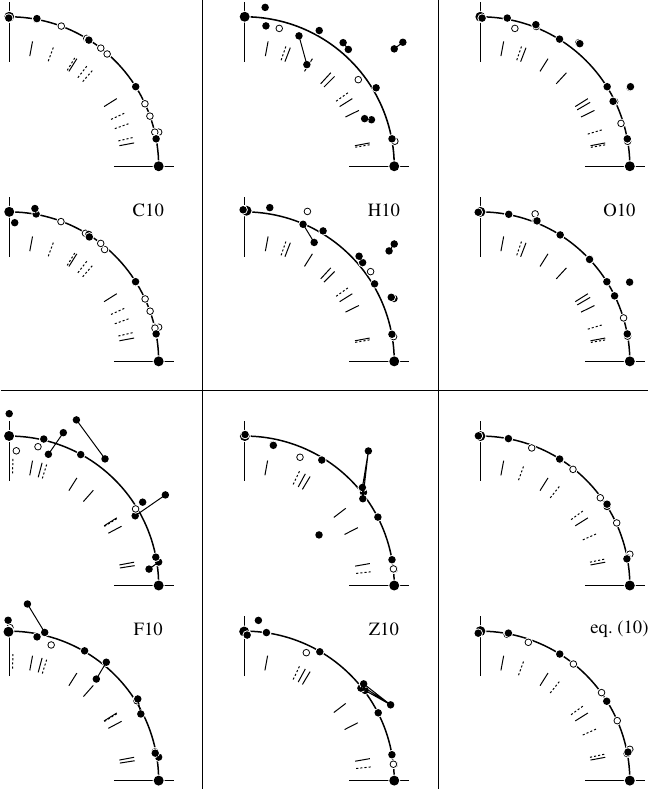}} \vspace{6pt} 
\caption{Application of one step $h = \pi / 2$ of C10, H10, O10, F10, 
Z10, and eq.~(\ref{the_method}) methods to systems of differential 
equations $\d x / \DT = -y$, $\d y / \DT = x$ (upper quarter-circles) and 
$\d x / \DT = -y / (x^2 + y^2)$, $\d y / \DT = x / (x^2 + y^2)$ (lower 
quarter-circles), with the initial condition $x(0) = 1$, $y(0) = 0$. The 
initial and final points, $(1 \SEP 0)$ and $\bigl( \tilde{x}(\pi/2) \SEP 
\tilde{y}(\pi/2) \bigr)$, respectively, are marked by large black dots. 
Smaller closed and open dots correspond to the intermediate positions 
$\X12S{i}\mkern0.25mu$, $1 \le i \le s$, with $\BI \ne 0$ and $\BI = 0$, 
respectively. Within node clusters, these are connected by thin solid 
lines. Nodes $\CI$, $1 \le i \le s$, are shown by radial ticks with solid 
($\BI \ne 0$) and dashed ($\BI = 0$) lines.} \label{circle} \end{figure}

\newpage\NI\begin{picture}(0,0)(0,0)\put(7,-534.5){\includegraphics[height=532.9pt]{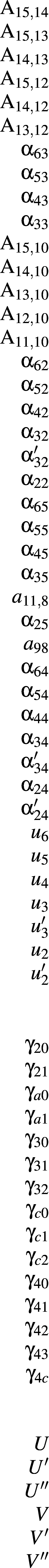}}\end{picture}\label{constants}\vspace{-1.5pt}
\CODE
          1        0        0        0       -1        0        0        0        2
         -107      353     -341      149     -821      1503    -509      843      1416
          815     -353      341     -149      821     -2211     509     -843      1416
         -1593179 -703001   523855   542527   7003159 -2315751 -2257457  672093   1458168
          325025  -242953  -352465  -38617    389855  -1114275 -170833  -1028451  1458168
          332873   157659  -28565   -83985   -1232169  693185   404715   59393    243028
          6       -2        0        1        12       3       -6        0        72
          30      -5       -6        1        30       21      -6       -15       144
          9        5        3       -3       -39       41       15       15       96
          11       1        1        1        3       -11      -3       -5        96
         -30363    37407   -20595    13773   -88581    226363  -12507    81737    55152
          3384    -11215    3141    -2166     11751   -51027    1596    -18702    18384
         -2167     7228    -1607     2060    -6999     35267   -1659     12133    27576
          117149  -60598    38062   -28367    199746  -411836   35073   -148213   55152
         -16257    10595   -5919     4243    -26211    67005   -6009     24579    13788
          467      244     -182     -88       1457    -580     -560      244      72
          551      318     -197     -126      2031    -734     -753      242      48
          1749     785     -525     -347      4857    -2475    -1653     639      864
         -17439   -8189     6135     3359    -52563    24135    19263   -7539     864
          637      289     -221     -119      1865    -879     -679      273      72
         -135     -97       51       37      -597      179      225     -65       24
         -1865     400      194     -442      5671     3900    -1474     294      72
          865     -520     -421      142     -2535    -846      1083     648      48
          4251     7555     4641     827      5007    -24357   -10047   -13155    864
          24411   -39763   -23727    6685    -149289   10101    75345    46851    864
         -1397     1219     715     -359      6259     1329    -2609    -1119     72
         -1002     373      249     -196      2934     1530    -939     -243      24
          55       28       20       11      -75      -175     -25      -65       2
         -12120    5558     4524    -2131     34092    16017   -12816   -5850     72
          18054   -15493   -6840     5831    -93786   -23529    35508    8961     144
         -4381     2147     1741    -769      12555    5507    -4839    -2359     96
          48165   -22829   -18405    8535    -137595  -62381    52215    24273    96
         -607      293      232     -110      1771     786     -671     -306      4
         -19287    12478    7254    -4733     75975    25299   -28668   -9450     72
          126      13       9        26      -270     -318      27      -69       18
          17771   -8607    -4789     4335    -64330   -28602    21794    4230     54
         -12118    5967     4931    -2052     33964    14784   -13307   -6807     36
          1153    -2902    -1712     323     -9149     2499     5122     3750     36
         -26894    22579    13153   -6818     118174   26544   -48563   -20445    36
          2881    -2231    -1307     727     -12243   -3213     4911     1917     12
          67950   -32635   -22419    13484   -212940  -92400    76203    27903    108
         -420     -198     -150     -80       574      1260     187      465      2

          45       10       6        4        31       54       10       18       18
         -1       -19      -1       -7       -49      -21      -19      -3        18
         -6       -3       -3       -1       -6       -21      -3       -9        12
          2        0        1        0        0        7        0        3        2
          1239     455      282      167      1134     2499     441      903      252
         -99      -125     -36      -41      -231     -375     -114     -141      36
          23       8        5        2       -7        42       5        24       12
         -14      -21      -16      -6       -21      -98      -21      -49       84
         -6        3        3        1        0        15       3        9        12
          1        0        0        0        0        0        0        0        1
          6048     2716     1725     1018     6622     15729    2479     5781     630
         -4480    -4039    -1498    -1384    -7840    -15708   -3472    -5943     630
          20      -2       -4       -2       -50      -24      -8        6        15
          20       9        8        3        14       63       8        27       15
          23       8        5        2       -7        42       5        24       15

          153566   36694    943      14484    178724   81235    53158   -28714    1
         -134733  -9751    -2187    -2447    -78969   -35907   -17925    14493    1
          1267700  0        0        0        0        0        0        0        1
          1854455 -323009   188367  -152875  -723177   1849113 -397521   765441   1
         -173229   32676   -14883    11874    142380  -149989   45426   -89567    1
          131201  -21770    7390    -4531    -154203   82173   -32718    73122    1
\end{Verbatim}

\NI\begin{picture}(0,0)(0,0)\put(0,-544.4){\includegraphics[height=541.3pt]{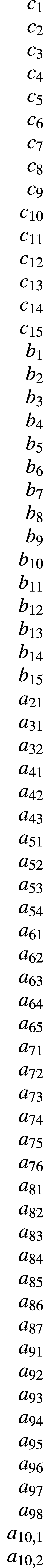}}\end{picture}\label{decimal}\vspace{-1.5pt}
\CODE
      +0.000000000000000000000000000000000000000000000000000000000000000000000000000000000000000000
      +0.133333333333333333333333333333333333333333333333333333333333333333333333333333333333333333
      +0.266666666666666666666666666666666666666666666666666666666666666666666666666666666666666667
      +0.400000000000000000000000000000000000000000000000000000000000000000000000000000000000000000
      +0.571428571428571428571428571428571428571428571428571428571428571428571428571428571428571429
      +0.778740761536291800442363524550407562438954342313508170658506769545868148577329534757889703
      +0.642615758240322548157075497020439535959501736363212695909875208263848965457099799090837864
      +0.882527661964732346425501486979669075182867844268052119663791177918527658519413257061748635
      +0.117472338035267653574498513020330924817132155731947880336208822081472341480586742938251365
      +0.117472338035267653574498513020330924817132155731947880336208822081472341480586742938251365
      +0.357384241759677451842924502979560464040498263636787304090124791736151034542900200909162136
      +0.357384241759677451842924502979560464040498263636787304090124791736151034542900200909162136
      +0.642615758240322548157075497020439535959501736363212695909875208263848965457099799090837864
      +0.882527661964732346425501486979669075182867844268052119663791177918527658519413257061748635
      +1.000000000000000000000000000000000000000000000000000000000000000000000000000000000000000000
      +0.033333333333333333333333333333333333333333333333333333333333333333333333333333333333333333
      +0.000000000000000000000000000000000000000000000000000000000000000000000000000000000000000000
      +0.000000000000000000000000000000000000000000000000000000000000000000000000000000000000000000
      +0.000000000000000000000000000000000000000000000000000000000000000000000000000000000000000000
      +0.000000000000000000000000000000000000000000000000000000000000000000000000000000000000000000
      +0.000000000000000000000000000000000000000000000000000000000000000000000000000000000000000000
      +0.000000000000000000000000000000000000000000000000000000000000000000000000000000000000000000
      +0.000000000000000000000000000000000000000000000000000000000000000000000000000000000000000000
      +0.135169627249231064398790288647151661598687390677589878805138875485701561628210909193328656
      +0.054067850899692425759516115458860664639474956271035951522055550194280624651284363677331462
      +0.215778257736022470617613537547175598111058915336253983819589520767421262523355528508005093
      +0.061650930781720705890746725013478742317445404381786852519882720219263217863815865288001455
      +0.277429188517743176508360262560654340428504319718040836339472240986684480387171393796006548
      +0.189237478148923490158306404106012326238162346948625830327194425679982186279495272870660119
      +0.033333333333333333333333333333333333333333333333333333333333333333333333333333333333333333
      +0.133333333333333333333333333333333333333333333333333333333333333333333333333333333333333333
      +0.000000000000000000000000000000000000000000000000000000000000000000000000000000000000000000
      +0.266666666666666666666666666666666666666666666666666666666666666666666666666666666666666667
      +0.100000000000000000000000000000000000000000000000000000000000000000000000000000000000000000
      +0.000000000000000000000000000000000000000000000000000000000000000000000000000000000000000000
      +0.300000000000000000000000000000000000000000000000000000000000000000000000000000000000000000
      +0.134110787172011661807580174927113702623906705539358600583090379008746355685131195335276968
      +0.000000000000000000000000000000000000000000000000000000000000000000000000000000000000000000
      +0.087463556851311953352769679300291545189504373177842565597667638483965014577259475218658892
      +0.349854227405247813411078717201166180758017492711370262390670553935860058309037900874635569
      +0.328658938997791240721282668696245180100949935428058492348520659513222104304668996831025694
      +0.000000000000000000000000000000000000000000000000000000000000000000000000000000000000000000
      -0.843123044711636120617820434267319755889422538822419939805710279613810068291349217990455336
      +1.230383064721160870667126440845906263177038685065274297872024392536130940503747853752249173
      +0.062821802528975809671774849275575875050388260642595320243671997110325172060261902165070172
      +0.130186990844533481397415277676876710451288921648802895948985973227427327921848463186124942
      +0.000000000000000000000000000000000000000000000000000000000000000000000000000000000000000000
      +0.000000000000000000000000000000000000000000000000000000000000000000000000000000000000000000
      +0.580738972917374107341609289970199905398640941037865280565410031214758908228344711848158528
      -0.132061017042928219570339015811796150446818327131853650470115559422384511796983327996330108
      +0.063750811521343178988389945185159070556390200808398169865594763244047241103889952052884502
      +0.116969144664032355608536137983019767284447494556651903864723553337776482921362147717802148
      +0.000000000000000000000000000000000000000000000000000000000000000000000000000000000000000000
      +0.000000000000000000000000000000000000000000000000000000000000000000000000000000000000000000
      +0.804791321585223079797049560718389631085997501079130322088670870267840258541270110417387198
      -1.032442407303393452179074692383412947850773751656102435103865475343877801453723942152851429
      +0.141237507577425159424271186625123531674337461239364982255644624276392633209837338226267645
      +0.851972095441445203774719294036549092988859139049007346558617605380396085300667602853143073
      +0.072301376815318053808169355285248577196274719038173036516820302048727944927561329493507003
      +0.000000000000000000000000000000000000000000000000000000000000000000000000000000000000000000
      +0.000000000000000000000000000000000000000000000000000000000000000000000000000000000000000000
      +0.278449353793857950736187552433349927872809528174204139729412569505988600192254101364671140
      -0.720355742529399485982635010432963793377590743017129914556141009831958711185177678815313348
      +0.071478309723007569077180467110101963232268708471224068127334024952865428728115593534703972
      +0.480546127100390254710165280839106235150089447220646221776588301884109281650681775624438024
      -0.064947086867906688774569132214511985256719504155169671257805366478260202832848378263755427
      -0.033604580075425596431640044865448312697754125457571033966446479388085421715398642048231498
      +0.000000000000000000000000000000000000000000000000000000000000000000000000000000000000000000
\end{Verbatim}

\newpage

\NI\begin{picture}(0,0)(0,0)\put(0,-544.4){\includegraphics[height=541.3pt]{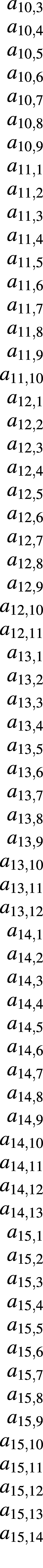}}\end{picture}\vspace{-1.5pt}
\CODE
      +0.000000000000000000000000000000000000000000000000000000000000000000000000000000000000000000
      -0.857669388668674924032812732229566068998144989807984562742043274869246367307599404275438590
      +2.054105395119576054328336908587945831582620175191578897288368999063154807831596816873316215
      -0.209909042661849185066865444237374444876978972036734389443755323786107747791309251353961678
      -1.333972235951611409051125524672704502068159831145103714187256702629992309144152509633176930
      +0.180806348384581339693160596923480075192005804935788561945038863986674416934017278251008686
      +0.317715841888671374135444753513998346683544094051974121442302739705074962673432455124735159
      +0.037810053861933727583877526864046322952305373764286754688878816133302931065788321936975530
      +0.000000000000000000000000000000000000000000000000000000000000000000000000000000000000000000
      +0.000000000000000000000000000000000000000000000000000000000000000000000000000000000000000000
      +0.045591558641579644208161372472470642588279099631100286855624026273895766568226149784507397
      +0.374858247294747357801811789348440678077393709386181864992340497763246020363342448831864027
      -0.018981904165238088665273073099056486944960468054380557714800392623159772511027805727133276
      -0.336038862637559379967733799098120208915305739878288300584069940705912757234040695695257751
      +0.042207716375559153022580951807545725160812786625465879746013555646386854001549779891730475
      +0.051900360983690832037790128319440699859433125187462257771821672390157423238366045436562647
      +0.160037071404964205821709606364793091262540376974959118334316556858234569050695956449913087
      +0.025378474265434522073341785109062596527860323276820540705890698859323818498838966304096456
      +0.000000000000000000000000000000000000000000000000000000000000000000000000000000000000000000
      +0.000000000000000000000000000000000000000000000000000000000000000000000000000000000000000000
      +0.362306327543474080870356990535692676714482207942118981416939928949791728367123473328004120
      -2.130022066259949711216588139786172085385079053283141571779776466074951544604054464037893691
      +0.167271594859437715452103295041620121357577378494618057193516457140141705837480748503078071
      +1.582750941864663213591839501273807297171503649177669091900056915258923612981224332245177665
      -0.204720185139081567587174987631558173178320360472369067149216918644283831827088288528863705
      +0.319265328144806137168538521838050367024648119039965645533219052127751799436983299054606914
      -0.199485256254696927298131093258193563570518501610697726802791703448058055735145490139989292
      +0.434639082735589988788638629857251227378344501071803353072286827567511801587537624180945598
      +0.021029555123492649231004219196592342935742304160151528452676708946550307466281435746860028
      +0.000000000000000000000000000000000000000000000000000000000000000000000000000000000000000000
      +0.000000000000000000000000000000000000000000000000000000000000000000000000000000000000000000
      -0.208531388344693562066394580214099278549298495126493308746134843353719483236004257692460587
      +0.326863498654767652665730419629359718148991455519984894344217473947412982044929978522731887
      -0.040291595063466341977981502499810508092167265608412827254081246436545949894568768120086225
      -0.062168570110085104851127006245165266954000899589268210982639466796612080783267242335332813
      +0.019401511635634607470015991594882193791752892731789828413457858797773885042894717414314220
      +0.033618934322149698701521634179531768990520248587976979368806692500394465464703080284812928
      +0.184057539066668989188409585677311173226149887887488899665117691435303709700110659619387592
      +0.126719042450063065291362679244566355420313628999369516069616211043815410052049067488162884
      +0.241917230505790894504534056457271037041497978800625396578838128179475719599971128162447951
      +0.059015439601144226763876106824299822651056803159364026683828389696455698859560833618210307
      +0.000000000000000000000000000000000000000000000000000000000000000000000000000000000000000000
      +0.000000000000000000000000000000000000000000000000000000000000000000000000000000000000000000
      +0.346753529400439458709721774328857076918185484739163574624551262985723353031416083991491991
      -0.543520439250941361637466822737163641756559230304227658933756185028034640683436881957640271
      +0.066998320513439386022349163118524383402374304344205128122905426297240312187142777357799209
      -0.209105926273890705326433037391968985907198847287185344061438195214198838012365442723920080
      +0.034196819215377393468216514309803928830340701888167708646927556886609080880789258681682710
      +0.656808691640117064523673744875271638362693883629496564722278370329924383275858838375005043
      -0.530023086491271629685194176217211641242245756479495833589514996293186125976018346545708951
      +0.682370713621575520444713157488350028626121649514411954149882119768342984225842629755650051
      -0.524316794649154201841271672103049178047360868487598331909107445957588634113426671093395740
      +0.843350394637897194983316734483955643345459719551750331207234874447240084844050177602574367
      -0.040440866328662678612786513864723972854343311900228087268921187067545862661799086441064095
      +0.000000000000000000000000000000000000000000000000000000000000000000000000000000000000000000
      +0.000000000000000000000000000000000000000000000000000000000000000000000000000000000000000000
      -0.936174412160392166960880644803504126455580435000519336040696957832662353800600021951728948
      +1.467410954959016212947954342116057777884286077126485056662059170859341898803083590224129382
      -0.180883849778987051699697942242124022922323870654198891654454984812850334688440411791417790
      +1.167589102274410575304328368479226265476383761610677844233764974596479007878520383510866872
      -0.280113604810353377498686142206136702751291429624401159362004219851994095769258196924953998
      -1.871668376348516420302371567652811494220817967866998659532896736116592408041514651376908974
      +2.241588116125051681252015179294157675223848655884912766574651463604345999471163122136357746
      -1.572572082691418559566908951416310106551328354276059431821310834756397117889598376197133068
      +2.151700399601058696148477230250852973515302452974155466886336399845840759500171462284652967
      -1.813340450902764402851969633637645578401662879405552525500029323314042291416407578145425752
      +0.666905070061557491840526275682961312057527301131726956823502234846076798614679764672625658
\end{Verbatim}

\end{document}